\documentclass[12pt]{article}
\usepackage{amsmath}
\usepackage{amsfonts}

\advance \textheight by \topskip
\advance \textheight by 1pt
\makeatother
\def\bea{\begin{eqnarray}}
\def\ena{\end{eqnarray}}

\def\lar{\longrightarrow}
\def\non{\nonumber}
\def\epn{\varepsilon_n}
\def\deg{\hbox{deg}}

\def\dim{\hbox{dim}}

\def\ord{\hbox{ord}}
\def\Res{\hbox{Res}}
\def\homega{\widehat{\omega}}
\def\Xaf{X^{{\rm aff}}}
\def\tX{\tilde{X}}
\def\tE{\tilde{E}}

\def\tB{\tilde{B}}
\def\bB{\bar{B}}
\def\tilp{\tilde{p}}
\def\tilq{\tilde{q}}
\def\tilp0{\tilde{p}_0}

\def\tip{p^{(i)}}
\def\tgamma{\tilde{\gamma}}
\def\tigamma{\tilde{\gamma}^{(i)}}
\def\igamma{\gamma^{(i)}}
\def\bgamma{\bar{\gamma}}
\def\tinfty{\tilde{\infty}}
\def\lar{\longrightarrow}

\def\hF{\widehat{F}}
\def\tilr{\tilde{r}}
\def\hr{\widehat{r}}

\newcommand{\bc}[2]{
\left(
\begin{array}{c}{#1}\\{#2}\end{array}
\right)}
\newcommand{\qbc}[2]{
\left[
\begin{array}{c}{#1}\\{#2}\end{array}
\right]}

\newcommand{\qed}{\hbox{\rule[-2pt]{3pt}{6pt}}}
\newtheorem{prop}{Proposition}
\newtheorem{theorem}{Theorem}
\newtheorem{defn}{Definition}
\newtheorem{lemma}{Lemma}

\title{
On Algebraic Expressions of Sigma Functions
for $(n,s)$ Curves
}

\author{
Atsushi Nakayashiki\thanks{
e-mail: 6vertex@math.kyushu-u.ac.jp}\\
Department of Mathematics, Kyushu University\\
}

\date{}
\begin{document}
\maketitle

\vskip2mm
\centerline{
Dedicated to Masaki Kashiwara on his sixtieth birthday
}
\vskip15mm

\begin{abstract}
An expression of the multivariate sigma function associated with a (n,s)-curve
is given in terms of algebraic integrals. As a corollary the first term of the
series expansion around the origin of the sigma function is directly proved 
to be the Schur function determined from the gap sequence at infinity.
\end{abstract}
\clearpage

\section{Introduction}
One of the prominent features of Weierstrass' elliptic sigma functions
is their algebraic nature directly related to the defining equation of the
elliptic curve. Klein \cite{K1,K2} extended the elliptic sigma functions to the
case of hyperelliptic curves from this point of view. 
Since they are defined, it had been one of the central problems to determine
the coefficients of the series expansion of the sigma functions.
This problem was studied mainly by making linear differential equations 
satisfied
by sigma functions \cite{W,Wil,Bo1,Bo2,Bo3,Bo4}, by making non-linear equations
\cite{B3} for genus two and by using algebraic expressions \cite{K1,K2,Burk}.
To determine a solution of linear differential equations it is necessary 
to specify an initial condition which requires separate 
consideration. Therefore the expansion was mainly studied for the sigma 
functions with non-singular even half periods as characteristics. 

Recently Klein's sigma function is further generalized to the case of more 
general plane algebraic curves called $(n,s)$-curves 
by Buchstaber, Enolski and Leykin \cite{BEL1,BEL2,BEL3,BEL4, BL}. 
They made an important observation that the first term, with respect to certain
degree introduced in the theory of soliton equations, of the series expansion
of the sigma function, which corresponds to the most singular characteristic,
is described by Schur function. Although such connection is expected 
from the theory of the KP-hierarchy \cite{Sat,DJKM},
a concrete description of the degeneration of the quasi-periodic 
solutions to singular curves has not been done before. 
In order to establish the connection to Schur functions 
Buchstaber et al.\cite{BEL2} have developed 
the rational theory of abelian integrals and characterized Schur functions 
by Riemann's vanishing theorem. 
Moreover Buchstaber and Leykin \cite{BL} have proposed a system of 
linear differential equations satisfied by sigma functions, which is independent of characteristics. Combining those results is expected to be effective
for the further study of sigma functions with singular characteristics.
Unfortunately a basic formula ( the formula (4) in \cite{BEL3}), on which
some of results of \cite{BL,BEL2} including the assertion related to 
Schur functions mentioned above depend, is not correct.

The purpose of this paper is to generalize Klein's algebraic formulas for the
hyperelliptic sigma functions to the case of $(n,s)$-curves and is to establish the relation with Schur functions directly, that is,  without using the results of \cite{BEL2,BEL3}. 
The sigma function, in this paper, signifies the sigma function with 
Riemann's constant as its characteristic \cite{BEL1}.

The heart of the algebraic formula for the extended sigma function is in the formula for the elliptic sigma function given by Klein \cite{K1}.
Therefore let us briefly explain it. 
Let $\sigma(u)$ and $\wp(u)$ be Weierstrass' sigma and elliptic functions
associated with the periods $2\omega_1$, $2\omega_2$.
Consider two variables $u_1$, $u_2$ and set 
$p_i=(x_i,y_i)=(\wp(u_i),\wp'(u_i))$.
They are points on the elliptic curve $y^2=4x^3-g_2x-g_3$.
By making use of the addition theorem for $\wp(u)$ the bilinear form

$$
\homega=\wp(u_2-u_1)du_1du_2
$$ 
can be written in an algebraic form as
\bea
&&
\homega=
\frac{2y_1y_2+4x_1x_2(x_1+x_2)-g_2(x_1+x_2)-2g_3}
{4y_1y_2(x_1-x_2)^2}dx_1dx_2.
\non
\ena
With this $\homega$ the elliptic sigma function is expressed as
\bea
&&
\sigma(u_2-u_1)=
\frac{x_1-x_2}{\sqrt{y_1y_2}}
\exp\left(\frac{1}{2}\int_{{\bar p}_1}^{{\bar p}_2}\int_{p_1}^{p_2}\homega\right),
\label{sigma-alg}
\ena
where ${\bar p}_i=(x_i,-y_i)$. 
What is remarkable for this formula is that the sigma function is expressed
by the algebraic functions $x_i,y_i$ and the integral of the algebraic
differential form $\homega$. 
This formula very clearly manifests the algebraic
structure of the sigma function.  
For example, prescribing degree $2i$ to $g_{i}$, one can deduce that 
the coefficients of the series expansion of $\sigma(u)$ at the origin 
become homogeneous polynomials of $g_2$ and $g_3$ directly from this formula 
without using differential equations.

For higher genus curves one needs to introduce $g$ variables in the sigma
function. Already in the case of genus one it is possible to introduce
arbitrary number of variables. In fact
the generalized addition formula due to
Frobenius and Stickelberger makes it possible to express the "$n$-point
function" in terms of the "$2$-point function":
\vskip7mm
$
\displaystyle{
\sigma\left(\sum_{i=1}^N(u_i-v_i)\right)
=
\frac
{
\prod_{i,j=1}^N\sigma(u_i-v_j)
\,\,\,\,
\det\left(\wp^{(i-1)}(u_j)\right)_{1\leq i,j\leq 2N}
}
{
\prod_{i<j}\sigma(u_i-u_j)\sigma(v_j-v_i)
\prod_{i,j=1}^N\left(\wp(v_j)-\wp(u_i)\right)
},
}
$
\vskip7mm
\noindent
where in the determinant we set $u_{N+j}=-v_j$, $1\leq j\leq N$. This
formula suggests how one should increase the number of variables in
the sigma function in general.

Consider the algebraic curve $X$ defined by
$$
y^n-x^s-\sum_{ni+sj<ns}\lambda_{ij}x^iy^j=0,
$$ 
with $n$ and $s$ being relatively coprime integers satisfying $2\leq n<s$.
We call it a $(n,s)$ curve. If it is non-singular its genus is
$g=1/2(n-1)(s-1)$.
The sigma function for $X$ is defined as the holomorphic function on 
${\mathbb C}^g$ which satisfies certain quasi-periodicity and 
normalization conditions (see (\ref{transformation}), (\ref{normalization})).
It can be considered as a holomorphic section of some line bundle on the
Jacobian $J(X)$ of $X$ or as a multi-valued holomorphic function on $J(X)$
whose multivalued property is specified by the quasi-periodicity. In turn
it can also be considered as a symmetric multi-valued holomorphic function
on $X^g$ through Abel-Jacobi map. More generally we construct a symmetric
multi-valued holomorphic function on $X^N$ with the required quasi-periodicity
properties for any $N\geq 1$.

The building block of the formula is the prime function which is a certain
modification of the prime form \cite{F}. It takes a similar form
to the right hand side of (\ref{sigma-alg}):
\bea
&&
\tE(p_1,p_2)=\frac{x(p_2)-x(p_1)}{\sqrt{f_y(p_1)f_y(p_2)}}
\exp\left(
\frac{1}{2}
\sum_{i=1}^{n-1}\int_{p_1^{(i)}}^{p_2^{(i)}}
\int_{p_1}^{p_2}
\homega
\right),
\non
\ena
for certain algebraic bilinear form $\homega$, where $\{p^{(0)},...,p^{(n-1)}\}$
is the inverse image of $p=p^{(0)}$ by the map $x: X \lar {\mathbb P}^1$,
$x:(x,y)\mapsto x$. It is skew symmetric and has the same transformation rule as that of the
sigma function when one of the argument goes round cycles of $X$.
Therefore this function can be considered as 
something like the sigma function 
restricted to the Abel-Jacobi image of $X\times X$ 
although the restriction of the sigma function itself vanishes identically.
Then the function on $X^N$ $(N\geq 3)$ is constructed in a form suggested
by Frobenius-Stickelberger's formula using the "$2$-point function"
$\tE(p_1,p_2)$. 
In this way the problem of constructing the sigma function reduces to 
finding certain meromorphic function on $X^N$. This problem is solved
for $(n,s)$-curve. 

To write explicitly the formula 
we need to describe a basis of meromorphic functions
on $X$ which are singular only at $\infty$. 
Prescribe degrees $n$ and $s$ to
$x$ and $y$, order the functions $x^iy^j$, $i\geq 0$, $0\leq j\leq n-1$ 
from lower degrees
and name them as $f_1=1$, $f_2$, $f_3$,...
Then the formula for the sigma function takes the form 
(Theorem \ref{sigma-2}):
\bea
&&
\sigma\left(\sum_{i=1}^N\int_{q_i}^{p_i}du\right)
=C_N\,\,M_N\,\,F_N,
\non
\ena
where $du$ is the vector of a basis of holomorphic one forms (\ref{1st-kind}),
$C_N$ is an explicit constant (\ref{cn}) and
\bea
M_N&=&
\frac
{
\prod_{i,j=1}^N\tE(p_i,q_j)
}
{
\prod_{i<j}
\left(
\tE(p_i,p_j)\tE(q_i,q_j)
\right)
\prod_{i,j=1}^N\left(x(p_i)-x(q_j)\right)
},
\non
\\
F_N&=&
\frac
{
D_N
}
{
\prod_{i<j}\left(x(q_i)-x(q_j)\right)^{n-2}
\prod_{k=1}^N\prod_{1\leq j\leq n-1}
\left(y(q_k^{(i)})-y(q_k^{(j)})\right)
},
\non
\\ 
D_N&=&\det\left(f_i(p_j)\right)_{1\leq i,j\leq nN},
\non
\ena
where we set 
$$
p_{N+(n-1)(k-1)+j}=q_k^{(j)},
\qquad 
1\leq k\leq N,
\quad  
1\leq j\leq n-1.
$$

In the case of hyperelliptic curves of genus $g$, 
that is, the case of $(n,s)=(2,2g+1)$, 
$F_N=D_N$, $q_k^{(1)}=(x,-y)$ for $q_k=(x,y)$ and 
the formula coincides with that given by Klein \cite{K2}.

It follows from this formula that, prescribing degrees $ns-ni-sj$ to 
$\lambda_{ij}$, the Taylor coefficients of the
sigma function become homogeneous polynomials of $\lambda_{ij}$
and the first term, with respect to certain degrees, of the expansion
of the sigma function is a Schur function corresponding to the partition
determined from the gap sequence at $\infty$ (Theorem \ref{main-cor}).

The plan of the present paper is as follows.
In section 2 necessary facts on Riemann surfaces and related objects on them
such as flat line bundles, prime form and normalized bilinear form
are reviewed. The meromorphic functions and
differentials on $(n,s)$-curves are studied in section 3. The Important object 
here is the algebraic bilinear form $\homega$. 
The existence of it is proved in section 3.3 and the relation
with the symplectic basis of the first cohomology group of a $(n,s)$-curve 
 is given in section 3.4. 
In section 4 the properties of Schur functions are reviewed.
The sigma function of a $(n,s)$-curve is defined and studied in section 5.
After giving the definition and an analytic expression of the sigma function
in section 5.1, an algebraic expression of the prime form is given 
in section 5.2.  In section 5.3 the prime function is introduced
and its properties are established using those of the prime form. 
The algebraic expressions of the sigma function are given in section 5.4.
Theorems \ref{sigma-1} and \ref{sigma-2} are main results of this paper. 
In section 5.5 the series expansion 
of the sigma function is studied and the proportionality constants 
in the proofs of main theorems are determined.
Examples of $(2,3)$ curve and more generally $(2,2g+1)$ curves
 are given in section 5.6 1nd 5.7. In section 6 some comments are given.

\section{Preliminaries}
\subsection{Riemann's Theta Function}
Let $\tau$ be a $g\times g$ symmetric matrix whose imaginary 
part is positive definite and
$a, b\in {\mathbb R}^g$. Riemann's theta function with
characteristics ${}^t(a,b)$ is defined by
\bea
&&
\theta\qbc{a}{b}(z)=\sum_{n\in {\mathbb Z}^g}
\exp\left(\pi i{}^t(n+a)\tau (n+a)+2\pi i(n+a)(z+b)\right).
\non
\ena
The theta function with zero characteristic is simply denoted by $\theta(z)$.
We list here some of the fundamental properties of Riemann's theta functions.
\noindent
(i)
\bea
&&
\theta\qbc{a}{b}(z+m_1+\tau m_2)
=
\non
\\
&&
\exp\left(
2\pi i({}^t a m_1-{}^t b m_2)-\pi i {}^t m_2\tau m_2-2\pi i {}^t m_2 z
\right)\theta\qbc{a}{b}(z),
\quad
m_1,m_2\in {\mathbb Z}.
\label{R-theta1}
\ena
\noindent
(ii)
\bea
&&
\theta\qbc{a}{b}(-z)=(-1)^{4{}^ta b}\theta\qbc{a}{b}(z),
\qquad
a,b\in \frac{1}{2}{\mathbb Z}^g
\label{R-theta2}
\ena
\noindent
(iii)
\bea
&&
\theta\qbc{a}{b}(z)=
\exp\left(
\pi i {}^t a\tau a+2\pi i {}^t a (z+b)
\right)
\theta(z+\tau a +b).
\label{R-theta3}
\ena

\subsection{Flat Line Bundle}
We briefly review some fundamental facts about Riemann surfaces 
and the description of flat line bundles on them.
We refer to \cite{F,G} for more details.

Let $X$ be a compact Riemann surface of genus $g$,
$\tX$ the universal cover of $X$ and $\pi$: $\tX\lar X$ be the projection.
We fix a marking of $X$. It means that we fix a base point $p_0$ on $X$,
 a base point $\tilp0$ on $\tX$ which lies 
over $p_0$ and a canonical basis $\{\alpha_i,\beta_j\}$ of $\pi_1(X,p_0)$, . 
Then the covering transformation group can be canonically identified with 
$\pi_1(X,p_0)$. 
For $k=0,1$, a holomorphic $k$-form on $X$ can be identified with that on $\tX$
which is $\pi_1(X,p_0)$-invariant.

Let $dv_j$, $1\leq j\leq g$ be the basis of holomorphic one forms normalized
as $\int_{\alpha_j}dv_j=\delta_{ij}$ and $\tau$ the
period matrix, $\tau=(\int_{\beta_j}dv_i)$. Set $dv={}^t(dv_1,...,dv_g)$.
The Jacobian variety $J(X)$ is defined by $J(X)={\mathbb C}^g/\tau {\mathbb Z}^g+{\mathbb Z}^g$. 

Let $S^kX=X^k/S_k$ be the $k$-th symmetric product of $X$. 
An element of it can be considered as a positive divisor on $X$ of degree $k$.
We denote by $I_k$ the Abel-Jacobi map with the base point $p_0$:
\bea
&&
I_k: S^kX\lar J(X),
\qquad
I_k(p_1+\cdots+p_k)=\sum_{i=1}^k \int_{p_0}^{p_i} dv.
\non
\ena
Then $J(X)$ can be identified with ${\rm Pic}^0(X)$ of linear equivalence
classes of divisors of degree zero by Abel-Jacobi map: for $A=\sum_{i=1}^d p_i$,$B=\sum_{i=1}^d q_i$,

\bea
I:{\rm Pic}^0(X)&\lar& J(X),
\non
\\
B-A &\mapsto& I(B-A)=I_d(B)-I_d(A).
\non
\ena
We sometimes use $I_k$ for the map $X^k\lar J(X)$.

A flat line bundle on $X$ is described by a representation
 $\chi:\pi_1(X,p_0)\lar {\mathbb C}^\ast$, 
where ${\mathbb C}^\ast$ is the multiplicative group of non-zero 
complex numbers.  
Namely a meromorphic section of the line bundle defined by $\chi$
is described by a meromorphic function $F$ on $\tX$ which satisfies
\bea
&&
F(\gamma \tilde{p})=\chi(\gamma)F(\tilde{p}).
\non
\ena
Since ${\mathbb C}^\ast$ is abelian, the image 
$\chi(\gamma)$ of $\gamma\in \pi_1(X,p_0)$ depends only on the image
of $\gamma$ in the homology group $H_1(X,{\mathbb Z})$, 
which we call the abelian image of $\gamma$.

Two representations $\chi_1$ and $\chi_2$ defines a holomorphically equivalent
line bundle if and only if 
\bea
&&
\chi_1(\gamma)\chi_2(\gamma)^{-1}=\exp\left(\int_\gamma\omega\right)
\non
\ena
for some holomorphic one form $\omega$ and any $\gamma\in \pi_1(X,p_0)$.

The Jacobian variety can also be identified with the set of 
holomorphic equivalence classes of flat line bundles on $X$.
The flat line bundle corresponding to the degree zero divisor $A-B$ with
$A$, $B$ positive divisors as before,
is described by
\bea
&&
\chi(\alpha_i)=1,
\qquad
\chi(\beta_i)=\exp\left(\int_A^Bdv_i\right),
\label{non-unitary}
\ena
where 
\bea
&&
\int_A^B dv=\sum_{i=1}^g \int_{p_i}^{q_i}dv,
\non
\ena
with the path from $p_i$ to $q_i$ being specified. Another choice of 
paths gives an equivalent line bundle.
We denote the equivalence class of this bundle by ${\cal L}_\alpha$, where 
$\alpha=\int_A^Bdv\in J(X)$. 

For $\alpha\in {\mathbb C}^g$ there exists a unique set of vectors
$\alpha',\alpha''\in {\mathbb R}^g$ such that
\bea
&&
\alpha=\tau \alpha'+\alpha''.
\non
\ena
The vector ${}^t(\alpha',\alpha'')$ is called the characteristic
of $\alpha$. We sometimes identify $\alpha$ with its characteristic.
Let ${}^t(\alpha',\alpha'')$ be the characteristic of 
$\int_A^Bdv\in {\mathbb C}^g$, where the integration paths are specified.
Then the function on $\tX$
\bea
&&
\frac{\theta(\int_{\tilde{p}_0}^{\tilde{p}}dv+\tau \alpha'+\alpha''+e)}
{\theta(\int_{\tilde{p}_0}^{\tilde{p}}dv+e)},
\non
\ena
is a meromorphic section of ${\cal L}_\alpha$ corresponding to $\chi$,
where $e\in {\mathbb C}^g$ is taken such that numerators and denominators
are not identically zero as a function of $\tilde{p}$.

There exists a unique unitary representation for each equivalence class of
line bundles.
The unitary representation for ${\cal L}_\alpha$
is given by
\bea
&&
\chi'(\alpha_j)=\exp(2\pi i \alpha'_j),
\qquad
\chi'(\beta_j)=\exp(-2\pi i \alpha''_j).
\label{unitary}
\ena

A Meromorphic sections of ${\cal L}_\alpha$ corresponding to $\chi'$
is given by 
\bea
&&
\frac{\theta[\alpha](\int_{\tilde{p}_0}^{\tilde{p}}dv+e)}
{\theta(\int_{\tilde{p}_0}^{\tilde{p}}dv+e)},
\non
\ena
where $e$ satisfies the same conditions as before.

\subsection{Prime Form}
Let $\delta_0$ be Riemann divisor for the choice $(p_0, \{\alpha_i,\beta_j\})$
and $L_0$ the corresponding holomorphic line bundle of degree $g-1$. 
For $\alpha\in J(X)$ set
$L_{\alpha}={\cal L}_{\alpha}\otimes L_0$.

There exists a non-singular odd half period $\alpha$ \cite{M2,F}.
By Riemann's theorem there is a unique divisor $p_1+\cdots+p_{g-1}$
such that
\bea
&&
\alpha=p_1+\cdots+p_{g-1}-\delta_0,
\non
\ena
in $J(X)$.
Considering the function $\theta[\alpha](\int_x^y dv)$ we see that 
the divisor of the holomorphic one form
\bea
&&
\sum_{i=1}^g \frac{\partial \theta[\alpha]}{\partial z_i}(0) dv_i(p)
\non
\ena
is $2\sum_{i=1}^{g-1} p_i$. Since $\alpha$ is non-singular, there is
a unique, up to constant, holomorphic section of $L_\alpha$ which 
vanishes on $p_1+\cdots+p_{g-1}$.
Thus there exists a holomorphic section $h_\alpha$ of $L_\alpha$
such that
\bea
&&
h_\alpha^2(p)= \sum_{i=1}^g \frac{\partial \theta[\alpha]}{\partial z_i}(0) 
dv_i(p).
\non
\ena
We use the same symbol $h_\alpha$ for the pull back of $h_\alpha$ 
to $\tX$. Then the prime form \cite{F,M2,B1} is defined as
\bea
&&
E(\tilde{p}_1,\tilde{p}_2)=
\frac{\theta[\alpha](\int_{\tilde{p}_1}^{\tilde{p}_2}dv)}
{h_\alpha(\tilde{p}_1)h_{\alpha}(\tilde{p}_2)},
\quad
\tilde{p}_1,\tilde{p}_2\in \tX.
\label{p-form}
\ena
By construction it vanishes to the first order at $\pi(\tilde{p}_1)=
\pi(\tilde{p}_2)$ and at no other divisors.
Let $\pi_j: X\times X\lar X$ be the projection to the $j$-th component
and $I_2:X\times X\lar J(X)$ Abel-Jacobi map.
Then $E(\tilde{p}_1,\tilde{p}_2)$ can be considered as a holomorphic section of the line bundle 
$\pi_1^\ast L_0^{-1}\otimes \pi_2^\ast L_0^{-1}\otimes I_2^\ast \Theta$
on $X\times X$, where $\Theta$ is the line bundle on $J(X)$ 
defined by the theta divisor $\Theta=\{\theta(z)=0\}$.
Notice that the prime form does not depend on the choice of $\alpha$.

We list some fundamental properties of the prime form.
\vskip2mm
\noindent
(i) $E(\tilde{p}_2,\tilde{p}_1)=-E(\tilde{p}_1,\tilde{p}_2)$.
\vskip2mm
\noindent
(ii) $E(\tilde{p}_1,\tilde{p}_2)=0$ $\Longleftrightarrow$
$\pi(\tilde{p}_1)=\pi(\tilde{p}_2)$.
\vskip2mm
\noindent
(iii) For $\tilde{p}\in \tX$ take a local coordinate $t$ around $\tilde{p}$. 
Then the expansion in $t(\tilde{p}_2)$ at $t(\tilde{p}_1)$ is of the form
\bea
&&
E(\tilde{p}_1,\tilde{p}_2)\sqrt{dt(\tilde{p}_1)dt(\tilde{p}_2)}=t(\tilde{p}_2)-t(\tilde{p}_1)
+O\left((t(\tilde{p}_2)-t(\tilde{p}_1))^3\right).
\non
\ena
\vskip2mm
\noindent
(iv) Consider the function
\bea
&&
F(\tilde{p})=\frac{E(\tilde{p},\tilde{p}_2)}{E(\tilde{p},\tilde{p}_1)},
\non
\ena
for $\tilde{p}_1, \tilde{p}_2\in \tX$.
If the abelian image of $\gamma\in \pi_1(X,p_0)$ is
 $\sum_{i=1}^g m_{1,i}\alpha_i+\sum_{i=1}^g m_{2,i}\beta_i$,
\bea
&&
F(\gamma \tilde{p})=
\exp\left(
-2\pi i {}^t m_2 \int_{\tilde{p}_1}^{\tilde{p}_2}dv
\right)F(\tilde{p}),
\non
\ena
where $m_i={}^t(m_{i,1},...,m_{i,g})$.

\subsection{Normalized Fundamental Form}\label{NSD}
Let $K_X$ be the canonical bundle of $X$. A section of
$\pi_1^\ast K_X\otimes \pi_2^\ast K_X$ is called a bilinear form
on $X\times X$ and a bilinear form $w(p_1,p_2)$ is called symmetric
if $w(p_2,p_1)=w(p_1,p_2)$. Since 
\bea
&&
H^0(X\times X,\pi_1^\ast K_X\otimes \pi_2^\ast K_X)\simeq
\pi_1^\ast H^0(X,K_X)\otimes \pi_2^\ast H^0(X,K_X),
\non
\ena
any holomorphic symmetric bilinear form can be written as
\bea
&&
\sum c_{ij} dv_i(p_1)dv_j(p_2),
\qquad
c_{ij}=c_{ji},
\label{SBF}
\ena
where $c_{ij}$'s are constants.

We denote by $\Delta$ the diagonal set of $X\times X$:
\bea
&&
\Delta=\{(p,p)\,|\, p\in X\,\}.
\non
\ena

\begin{defn}\label{NFF}
A meromorphic bilinear form $\omega(p_1,p_2)$ on $X\times X$ 
is called a normalized fundamental
form if the following conditions are satisfied.
\vskip2mm
\noindent
(i) $\omega(p_1,p_2)$ is holomorphic except $\Delta$ where it has a double pole. For $p\in X$ take a local 
coordinate $t$ around $p$. Then the expansion in $t(p_1)$ at $t(p_2)$ is 
of the form
\bea
&&
\omega(p_1,p_2)=
\left(\frac{1}{(t(p_1)-t(p_2))^2}+\hbox{regular}\right)dt(p_1)dt(p_2).
\label{omega-exp}
\ena
\vskip2mm
\noindent
(ii) $\displaystyle{\int_{\alpha_j}\omega=0}$, where the integration
is with respect to any one of the variables.
\vskip5mm
\end{defn}

Normalized fundamental form exists and unique.
It can be expressed explicitly using the prime form as \cite{F}
\bea
&&
\omega(p_1,p_2)=
d_{\tilde{p}_1}d_{\tilde{p}_2}
\log\,E(\tilde{p}_1,\tilde{p}_2),
\label{normalized-omega}
\ena
where $p_i=\pi(\tilde{p}_i)$.
Integrating this formula we get

\begin{prop}\label{exp-omega}\cite{F,G}
For $\tilde{a},\tilde{b},\tilde{c},\tilde{d}\in \tX$,
\bea
&&
\exp\left(\int_{\tilde{a}}^{\tilde{b}}\int_{\tilde{c}}^{\tilde{d}}\omega\right)
=
\frac{E(\tilde{b},\tilde{d})E(\tilde{a},\tilde{c})}{E(\tilde{a},\tilde{d})
E(\tilde{b},\tilde{c})}.
\non
\ena

\end{prop}
\vskip5mm

\section{$(n,s)$ curve}
\subsection{Definition}
For relatively coprime integers $n$ and $s$ satisfying
$s>n\geq 2$ consider the polynomial \cite{BEL3}
\bea
&&
f(x,y):=y^n-x^s-\sum_{in+js<ns} \lambda_{ij}x^iy^j.
\label{ns-curve}
\ena
Let $\Xaf$ be the plane algebraic curve defined by $f(x,y)=0$.
We assume that $\Xaf$ is non-singular.
Denote $X$ the corresponding compact Riemann surface
which can be considered as $\Xaf$ completed by one point $\infty$. 
The point $\infty$ becomes a ramification
point with the ramification index $n$. The genus of $X$ becomes
$g=1/2(n-1)(s-1)$. 
Hereafter we take $\infty$ as a base point and fix a marking of $X$, 
$(\infty,\tinfty,\{\alpha_i,\beta_i\})$.
 
A basis of holomorphic one form on $X$ is given by
\bea
&&
du_{i}=-\frac{x^{a_i-1}y^{n-1-b_i}dx}{f_y},
\label{1st-kind}
\ena
where $\{(a_i,b_i)\}$ is the set of non-negative integers satisfying 
\bea
&&
1\leq b\leq n-1,
1\leq a\leq [\frac{sb-1}{n}],
\non
\ena
and ordered as $-na_1+sb_1<\cdots<-na_g+sb_g$ \cite{BEL2}.
This order is specified in such a way that the order of zeros at
$\infty$ is increasing.
\vskip5mm

\noindent
{\large \bf Example}
\hskip2mm
$\displaystyle{du_g=-\frac{dx}{f_y}}$,
\hskip2mm
$\displaystyle{du_{g-1}=-\frac{xdx}{f_y}}$.

\subsection{Meromorphic Functions on $X$}
The space of meromorphic functions on $X$ which are holomorphic on 
$X\backslash\{\infty\}$ coincides with the space of polynomials of 
$x$ and $y$.
We describe a basis of this space.
Let $w_1<\cdots<w_g$ be the gap sequence at $\infty$.
It means that there is no meromorphic function on $X$ which has poles only 
at $\infty$ of order $w_i$.
Then 
\begin{lemma}\label{gap-seq}{\rm \cite{BEL2}}
\noindent
(i) $w_1=1$ and $w_{g}=2g-1$.
\vskip2mm
\noindent
(ii) Let $0=w_1^\ast<\cdots<w_g^\ast$ be integers such that
$\{w_i^\ast,w_i|i=1,...g\}=\{0,1,...,2g-1\}$. Then
$(2g-1-w_1^\ast,...,2g-1-w_g^\ast)=(w_g,...,w_1)$.
\end{lemma}

Notice that $\{w_i^\ast\}$ are non-gaps between $0$ and $2g-2$.

A local parameter $t$ around $\infty$ can be taken in such a way
that
\bea
&&
x=\frac{1}{t^n},
\qquad
y=\frac{1}{t^s}\left(1+O(t)\right).
\label{par-inf}
\ena
In particular $x$ and $y$ have poles at $\infty$ of order $n$ and $s$
respectively. For a meromorphic function $h$ on $X$ we denote
by ${\rm ord}_\infty \, h$ the order of poles at $\infty$.
Then
\bea
&&
{\rm ord}_\infty\, x^iy^j=ni+sj.
\non
\ena

Let $L(k\infty)$ be the vector space of meromorhic functions on $X$
which are holomorphic on $X\backslash\{\infty\}$ and have poles at
$\infty$ of order at most $k$. 
Set $L(\ast\infty)=\cup_{k=0}^\infty L(k\infty)$,
which is the space of meromorphic functions on $X$ holomorphic outside
$\infty$. A basis of $L(\ast \infty)$ is given by
\bea
&&
x^iy^j,
\qquad
i\geq 0,
\qquad
0\leq j\leq n-1.
\label{mon-basis}
\ena
There are exactly $g$-gaps in the set $\{{\rm ord}_\infty x^iy^j\}$:
\bea
&&
{\mathbb Z}_{\geq 0}
\backslash 
\{\,ni+sj\,|i\geq 0,\,0\leq j\leq n-1\,\}
=\{w_1<\cdots<w_g\}.
\non
\ena
Let $f_i$ be the monomial basis (\ref{mon-basis})
such that
\bea
&&
0={\rm ord}_\infty f_1<{\rm ord}_\infty f_2<{\rm ord}_\infty f_3<\cdots.
\non
\ena
In particular $f_1=1$, $f_2=x$.
By Riemann-Roch theorem
\bea
&&
\dim\, L((N+g-1)\infty)=N
\qquad
\hbox{for}
\quad
N\geq g.
\non
\ena
Explicitly, using the local coordinate $t$, we have
\bea
&&
f_i=
\left\{
\begin{array}{ll}
\displaystyle{\frac{1}{t^{w_i^\ast}}(1+O(t))} 
& 
\qquad 1\leq i\leq g
\\
[0.2cm]
\displaystyle{\frac{1}{t^{g-1+i}}(1+O(t))} 
& 
\qquad g+1\leq i
\end{array}
\right.
\label{fi-exp}
\ena
Notice that $\ord_{\infty}\,f_g=2g-2$.

\vskip5mm
\noindent
{\bf Example} $(n,s)=(2,2g+1)$:
\vskip2mm
$(w_1,...,w_g)=(1,3,...,2g-1)$, 
\hskip3mm 
$(w_1^\ast,...,w_g^\ast)=(0,2,4,...,2g-2)$,\par
\vskip2mm
$(f_1,f_2,...)=(1,x,...,x^{g},y,x^{g+1},xy,x^{g+2},x^2y,...)$.
\vskip3mm

\noindent
$(n,s)=(3,4)$, $g=3$:
\vskip2mm
$(w_1,w_2,w_3)=(1,2,5)$, 
\hskip3mm
$(w_1^\ast,w_2^\ast,w_3^\ast)=(0,3,4)$,
\vskip2mm
$(f_1,f_2,...)=(1,x,y,x^2,xy,y^2,x^3,x^2y,xy^2,...)$.
\vskip2mm

\noindent
$(n,s)=(3,5)$, $g=4$:\par
\vskip2mm
$(w_1,w_2,w_3,w_4)=(1,2,4,7)$, \hskip3mm
$(w_1^\ast,w_2^\ast,w_3^\ast,w_4^\ast)=(0,3,5,6)$,
\vskip2mm
$(f_1,f_2,...)=(1,x,y,x^2,xy,x^3,y^2,x^2y,x^4,xy^2,x^3y,...)$.
\vskip2mm

\noindent
$(n,s)=(3,7)$, $g=6$:\par
\vskip2mm
$(w_1,...w_6)=(1,2,4,5,8,10)$, 
\hskip3mm
$(w_1^\ast,...,w_6^\ast)=(0,3,6,7,9,10)$,
 \vskip2mm
$(f_1,f_2,...)=(1,x,x^2,y,x^3,xy,x^4,x^2y,y^2,...)$.
\vskip2mm

\noindent
$(n,s)=(4,5)$, $g=6$:\par
\vskip2mm
$(w_1,...w_6)=(1,2,3,6,7,11)$,
\hskip3mm 
$(w_1^\ast,...,w_6^\ast)=(0,4,5,8,9,10)$,
 \vskip2mm
$(f_1,f_2,...)=(1,x,y,x^2,xy,y^2,x^3,x^2y,xy^2,...)$.
\vskip5mm

In terms of $f_i$ the holomorphic one form $du_i$ is simply
described as 
\bea
&&
du_i=-\frac{f_{g+1-i}}{f_y}dx,
\qquad
1\leq i\leq g.
\non
\ena
By Lemma \ref{gap-seq} (ii) we have, around $\infty$,
\bea
&&
du_i=\left( t^{w_i-1}+O(t^{w_i})\right)dt.
\label{dui-exp}
\ena

\subsection{Algebraic Fundamental Form}
A meromorphic symmetric bilinear form which satisfies the condition (i)
of Definition \ref{NFF} can explicitly be constructed in terms of algebraic
functions.
Such algebraic form plays a central
role in the construction of the sigma function.

Let $p_i=(x_i,y_i)$, $i=1,2$ and 
\bea
&&
\Omega(p_1,p_2)=
\frac{\sum_{i=0}^{n-1}
y_1^i[\frac{f(z,w)}{w^{i+1}}]_{+}\vert_{(z,w)=(x_2,y_2)}}
{(x_1-x_2)f_y(p_1)}dx_1,
\non
\ena
where 
\bea
&&
[\sum_{n\in{\mathbb Z}}a_nw^n]_{+}=\sum_{n\geq 0}a_nw^n.
\non
\ena
Consider 
\bea
&&
\homega(p_1,p_2)=d_{p_2}\Omega(p_1,p_2)+
\sum c_{i_1j_1;i_2j_2}
\frac{x_1^{i_1}y_1^{j_1}}{f_y(p_1)}
\frac{x_2^{i_2}y_2^{j_2}}{f_y(p_2)}dx_1dx_2,
\label{h-omega}
\ena
where $(i_1,j_1)$ runs over $(a_i-1,n-1-b_i)$, $1\leq i\leq g$ 
and $i_2\geq 0$, $0\leq j_2\leq n-1$, $c_{i_1j_1;i_2j_2}$ 's 
are constants.
Assign degrees as
\bea
&&
\deg \,\lambda_{ij}=ns-ni-sj,
\quad
\deg\, x=\deg \,dx=n,
\quad
\deg\, y=s.
\non
\ena

\begin{prop}\label{cijkl}
\noindent
(i) If $c_{ij;kl}$ is taken such that $\homega(p_1,p_2)=\homega(p_2,p_1)$ then
$\homega$ defined by (\ref{h-omega}) satisfies the conditions (i) of Definition
\ref{NFF}.
\vskip2mm
\noindent
(ii) There exists a set of $c_{i_1j_1;i_2j_2}$ 
such that $\homega(p_1,p_2)=\homega(p_2,p_1)$,
non-zero  $c_{i_1j_1;i_2j_2}$ is a homogeneous polynomial of 
$\{\lambda_{kl}\}$ of degree $2sn-n(i_1+i_2+2)-s(j_1+j_2+2)$ and 
$c_{i_1j_1;i_2j_2}=0$ if $2ns-n(i_1+i_2+2)-s(j_1+j_2+2)<0$.
\end{prop}

Notice that, if we take $c_{i_1j_1;i_2j_2}$ as in (ii) in the proposition,
then $\homega$ becomes homogeneous of degree $0$.
For the proof of the proposition we need several lemmas.

Let $B$ be the set of branch points for the map
$x:X\lar {\mathbb P}^1$, $(x,y)\mapsto x$. 
For $p\in X$ set 
$x^{-1}(x(p))=\{p^{(0)},...,p^{(n-1)}\}$ with $p=p^{(0)}$, where the same $p^{(i)}$ is listed according to its multiplicity.

\begin{lemma}\label{lem-Omega1}
The one form $\Omega(p_1,p_2)$ is holomorphic except
$\Delta\cup \{(p^{(i)},p)|p\in B, i\neq 0\}\cup X\times \{\infty\}
\cup \{\infty\}\times X$.
\end{lemma}
\vskip2mm
\noindent
{\it Proof.}
It is sufficient to prove that $\Omega$ does not have a pole at $p_1=p_2^{(i)}$,$i\neq 0$, for $p_2\notin B$. 
Let $p_2^{(i)}=(x_2,y_2^{(i)})$ and
\bea
&&
f(x,y)=\sum_{j=0}^n f_j(x) y^j,
\quad
f_n=1.
\label{f+}
\ena
Then, for $i\neq 0$,
\bea
\sum_{k=0}^{n-1}(y_2^{(i)})^k
[\frac{f(z,w)}{w^{k+1}}]_{+}|_{(z,w)=(x_2,y_2)}
&=&
\sum_{k=0}^{n-1}(y_2^{(i)})^k\sum_{j\geq k+1}f_j(x_2)y_2^{j-k-1}
\non
\\
&=&
\sum_{j=1}^n f_j(x_2)\frac{(y_2^{(i)})^j-y_2^j}{y_2^{(i)}-y_2}
\non
\\
&=&
\frac{f(x_2,y_2)-f(x_2,y_2^{(i)})}{y_2^{(i)}-y_2}
\non
\\
&=&0,
\ena
where we use $y_2\neq y_2^{(i)}$, $i\neq 0$ which follows from
the assumption $p_2\notin B$.
Thus $\Omega$ is holomorphic at $p_1=p_2^{(i)}$, $i\neq 0$ as desired.
\qed

\begin{lemma}\label{lem-Omega2}
Let $p\notin B$, $t$ a local coordinate around $p$ and $t_i=t(p_i)$.
Then the expansion of $\Omega$ in $t_2$ at $t_1$ is of the form
\bea
&&
\Omega(p_1,p_2)=
\left(
\frac{-1}{t_2-t_1}+O\left((t_2-t_1)^0\right)
\right)dt_1.
\label{exp-Omega}
\ena
\end{lemma}
\vskip2mm
\noindent
{\it Proof.}
Since $p\notin B$ one can take $x$ as a local coordinate around $p$.
Therefore it is sufficient to prove
\bea
&&
\sum_{k=0}^{n-1}y_1^k
[\frac{f(z,w)}{w^{k+1}}]_{+}|_{(z,w)=(x_1,y_1)}=f_y(x_1,y_1).
\label{f+=fy}
\ena
Let us write $f(x,y)$ as in (\ref{f+}). Then the left hand side
of (\ref{f+=fy}) is equal to
\bea
&&
\sum_{k=0}^{n-1}y_1^k\sum_{j=k+1}^n f_j(x_1)y_1^{j-k-1}
=\sum_{j=1}^n jf_j(x_1)y_1^{j-1}=f_y(x_1,y_1).
\non
\ena
\qed
\vskip5mm

\begin{lemma}\label{lem-Omega3}
The meromorphic bilinear form $d_{p_2}\Omega(p_1,p_2)$ is holomorphic
except 
$\Delta \cup \{(p^{(i)},p)|p\in B, i\neq 0\}
\cup X\times \{\infty\}$.
\end{lemma}
\vskip2mm
\noindent
{\it Proof.}
Due to Lemma \ref{lem-Omega1} it is sufficient to prove
that $d_{p_2}\Omega(p_1,p_2)$ is holomorphic at 
$(\infty,p_2)$, $p_2\neq \infty$.

Let $t$ be the local coordinate around $\infty$ such that 
$x=1/t^n$, $y=(1/t^s)(1+O(t))$ and $t_i=t(p_i)$.
Then at $(\infty,p_2)$, $p_2\neq \infty$, the expansion
of $\Omega$ in $t_1$ takes the form
\bea
&&
\Omega=-\frac{dt_1}{t_1} \left(1+O(t_1)\right).
\non
\ena
Thus $d_{p_2}\Omega$ is holomorphic at $(\infty,p_2)$, $p_2\neq \infty$.
\qed
\vskip5mm

\begin{lemma}\label{lem-Omega4}
There exist second kind differentials $d\hr_i$, $1\leq i\leq g$
which are holomorphic outside $\{\infty\}$ and satisfy the
equation
\bea
&&
\omega(p_1,p_2)-d_{p_2}\Omega(p_1,p_2)=\sum_{i=1}^g du_i(p_1)d\hr_i(p_2).
\label{o-do}
\ena
\end{lemma}
\vskip5mm
\noindent
{\it Proof.}
Let us set
\bea
&&
\omega_1(p_1,p_2)=\omega(p_1,p_2)-d_{p_2}\Omega(p_1,p_2).
\non
\ena
By Lemma \ref{lem-Omega1}, \ref{lem-Omega2}, \ref{lem-Omega3} and 
(\ref{omega-exp}), the singularities of $\omega_1$
are contained in $B_2 \cup X\times \{\infty\}$, where 
$B_2=\{(b^{(i)},b)|b\in B\backslash\{\infty\}, 0\leq i\leq n-1\}$. 
Since $B_2$ is a finite set and $B_2\cap (X\times \{\infty\})=\phi$, $\omega_1$ is holomorphic except $X\times \{\infty\}$.
Thus one can write, for $p_2\neq \infty$, 
\bea
&&
\omega_1(p_1,p_2)=\sum_{i=1}^g du_i(p_1) d\tilr_i(p_2),
\non
\ena
for some one forms $d\tilr_i$. Let us describe $d\tilr_i$ 
more neatly in terms of $\omega(p_1,p_2)$. 
To this end let us take $q_1,...,q_g\in X\backslash B$ such that
$\sum_{j=1}^g q_i$ is a general divisor and $q_j$'s are in some small 
neighborhood of $\infty$. 
Take the local coordinate $t$ around $\infty$ 
such that $x=1/t^n$, $y=(1/t^s)(1+O(t))$ and write
\bea
du_i(p)&=&h_i(t)dt,
\non
\\
\omega_1(p_1,p_2)&=&K_1(t(p_1),p_2)dt(p_1).
\non
\ena
Then we have a set of linear equations
\bea
&&
\sum_{i=1}^g h_i(t(q_j))d\tilr_i(p_2)=K_1(t(q_j),p_2).
\non
\ena
Since $\sum_{j=1}^g q_j$ is a general divisor, $\det(h_i(t(q_j)))\neq 0$.
Thus $d\tilr_i$ can be expressed as
\bea
&&
d\tilr_i(p_2)=\sum c_{ij}K_1(t(q_j),p_2)
\non
\ena
for some constants, with respect to $p_2$, $c_{ij}$.

Notice that $K_1(t(q_j),p_2)$ is a second kind differential whose
only singularity is $\infty$. 
In fact the coefficient of $\omega(q_j,p_2)$ of $dt(p_1)$ is a second kind differential due to the property (i) of $\omega$ and $d_{p_2}\Omega(q_j,p_2)$ is
obviously of second kind. Thus $K_1(t(q_j),p_2)$ is a second kind differential.
As already proved the only singularity of $K_1(t(q_j),p_2)$ is $p_2=\infty$.

Let us set 
\bea
&&
d\hr_i(p)=\sum_{i=1}c_{ij}K_1(t(q_j),p),
\non
\ena
which is a second kind differential singular only at $\infty$, and set
\bea
&&
\omega_2(p_1,p_2)=\omega_1(p_1,p_2)-\sum_{i=1}^g du_i(p_1)d\hr_i(p_2).
\non
\ena
Then $\omega_2=0$ on $X\times (X\backslash\{\infty\})$.
Thus $\omega_2=0$ on $X\times X$.
Consequently 
\bea
&&
\omega(p_1,p_2)-d_{p_2}\Omega(p_1,p_2)=
\sum_{i=1}^g du_i(p_1) d\hr_i(p_2),
\non
\ena
which proves the lemma.
\qed
\vskip10mm

\noindent
{\it Proof of Proposition \ref{cijkl}}.
\vskip2mm
\noindent
(ii) Let us write
\bea
&&
d_{p_2}\Omega(p_1,p_2)=
\frac{\sum_{j_1,j_2\leq n-1}a_{i_1j_1;i_2j_2}x_1^{i_1}y_1^{j_1}x_2^{i_2}y_2^{j_2}}{(x_1-x_2)^2f_y(p_1)f_y(p_2)}dx_1dx_2.
\non
\ena
It can be easily verified that $a_{i_1j_1;i_2j_2}\in 
{\mathbb Z}[\{\lambda_{kl}\}]$ and $a_{i_1j_1;i_2j_2}$ is homogeneous
of degree $2(n-1)s-n(i_1+i_2)-s(j_1+j_2)$.

On the other hand
\bea
&&
\sum c_{i_1j_1;i_2j_2} 
\frac{x_1^{i_1}y_1^{j_1}x_2^{i_2}y_2^{j_2}}{f_y(p_1)f_y(p_2)}
=
\frac{\sum ( c_{i_1-2,j_1;i_2j_2}-2c_{i_1-1,j_1;i_2-1,j_2}+c_{i_1,j_1;i_2-2,j_2})x_1^{i_1}y_1^{j_1}x_2^{i_2}y_2^{j_2}}{(x_1-x_2)^2f_y(p_1)f_y(p_2)}.
\non
\ena
Thus $\homega(p_1,p_2)=\homega(p_2,p_1)$ is equivalent to
\bea
&&
c_{i_1-2,j_1;i_2j_2}-2c_{i_1-1,j_1;i_2-1,j_2}+c_{i_1,j_1;i_2-2,j_2}
-
c_{i_2-2,j_2;i_1j_1}+2c_{i_2-1,j_2;i_1-1,j_1}-c_{i_2,j_2;i_1-2,j_1}
\non
\\
&&
=a_{i_2,j_2;i_1j_1}-a_{i_1j_1;i_2j_2}.
\label{lineq1}
\ena
This is a system of linear equations for $c_{i_1j_1;i_2j_2}$ whose
coefficient matrix has integers as components.
\vskip5mm

\begin{lemma}\label{lem-Omega5}
Any meromorphic differential on $X$ which is singular only at $\infty$
is a linear combination of $(x^iy^j/f_y)dx$, $i\geq 0$, $0\leq j\leq n-1$.
\end{lemma}
\vskip2mm
{\it Proof.}
Let $\eta$ be a meromorphic differential which has a pole at $\infty$ of order
$k$ and is holomorphic on $X\backslash \{\infty\}$.
Consider the meromorphic function $\eta f_y/dx$. It has a pole only at $\infty$
of order $k-(2g-2)$ since the zero divisor of $dx/f_y$ is $(2g-2)\infty$.
Any meromorphic function holomorphic except $\infty$ is
a polynomial of $x$ and $y$. Thus $\eta\in {\mathbb C}[x,y]dx/f_y$.
\qed
\vskip5mm

By Lemma \ref{lem-Omega4} and \ref{lem-Omega5}, the system of linear equations
(\ref{lineq1}) has a solution, that is, it consists of compatible equations. 
Moreover it has a solution
such that each $c_{i_1j_1;i_2j_2}$ is a linear combination
of $a_{i_1'j_1';i_2'j_2'}$ satisfying $i_1+i_2=i_1'+i_2'+2$, 
$j_1+j_2=j_1'+j_2'$. In particular one can set $c_{i_1j_1;i_2j_2}=0$ 
if $2ns-n(i_1+i_2+2)-s(j_1+j_2+2)<0$ and has
\bea
&&
\deg\, c_{i_1j_1;i_2j_2}=2ns-n(i_1+i_2+2)-s(j_1+j_2+2),
\non
\ena
if $c_{i_1j_1;i_2j_2}$ is non-zero.
\qed

\vskip2mm
\noindent
(i) It is sufficient to prove the property (i) for $\omega$.
By Lemma \ref{lem-Omega4} $d_{p_2}\Omega(p_1,p_2)$ is holomorphic
except 
$\{(p,p)|p\in X\}\cup X\times \{\infty\}$ and so is $\homega$. Since 
$\homega(p_1,p_2)=\homega(p_2,p_1)$, $\homega$ does not have a pole at 
$p_2=\infty$ and therefore is holomorphic except $\{(p,p)|p\in X\}$ where 
it has a double pole. 

Let us prove that the expansion of
$\homega$ at the diagonal has the required form.
Set
\bea
&&
dr_i=-\sum c_{a_i-1,n-1-b_i;kl}\frac{x^ky^l}{f_y}dx.
\label{dr}
\ena
Then
\bea
&&
\homega-\omega=\sum_{i=1}^g du_i(p_1)(dr_i(p_2)-d\hr_i(p_2)).
\label{ho-o}
\ena
Both hand sides of (\ref{ho-o}) are meromorphic on $X\times X$.
The singularities of the left hand side are contained in
$\{(p,p)|p\in X\}$ and those of the right hand side are contained
in $X\times \{\infty\}$. Thus the possible singularity of $\homega-\omega$
is $\{\infty\}\times \{\infty\}$. Therefore $\homega-\omega$  and $dr_i-d\hr_i$
are holomorphic on $X\times X$ and $X$ respectively. Then the required expansion of $\homega$ at $(p,p)$, $p\in X$ follows from (\ref{omega-exp}).
\qed
\vskip5mm

Take one set of $c_{i_1j_1;i_2j_2}$ satisfying (ii) of the proposition 
and define $dr_i$ by (\ref{dr}).
Notice that $dr_i$ is a second kind differential.
In fact $dr_i=d\hr_i$ modulo holomorphic one form as is just proved 
and $d\hr_i$ is a second kind differential by Lemma \ref{lem-Omega4}.
We have
\bea
&&
\homega(p_1,p_2)=d_{p_2}\Omega(p_1,p_2)+\sum_{i=1}^gdu_i(p_1)dr_i(p_2).
\non
\ena
Define period matrices 
$\omega_1,\omega_2,\eta_1,\eta_2$ by
\bea
&&
2\omega_1=\left(\int_{\alpha_j}du_{i}\right),
\quad
2\omega_2=\left(\int_{\beta_j}du_{i}\right),
\quad
-2\eta_1=\left(\int_{\alpha_j}dr_{i}\right),
\quad
-2\eta_2=\left(\int_{\beta_j}dr_{i}\right).
\non
\ena
Notice that $\omega_1$ is invertible due to Riemann's inequality.
We set $\tau=\omega_1^{-1}\omega_2$. It is symmetric and satisfies 
$\rm{Im}\,\tau>0$.

\subsection{Relation between $\omega$ and $\homega$}
We give the relation between $\omega$ and $\homega$ using the period matrices.

\begin{lemma}\label{lem-oo1}
Let $\omega_1$ and $\omega_2$ be meromorphic symmetric bilinear form
satisfying the condition (i) of Definition \ref{NFF}. Then
\bea
&&
\omega_1-\omega_2=\sum_{i,j=1}^g c_{ij} du_i(p_1)du_j(p_2),
\label{o1-o2}
\ena
for some constants $c_{ij}$ such that $c_{ij}=c_{ji}$.
\end{lemma}
\vskip2mm
\noindent
{\it Proof.}
The left hand side of (\ref{o1-o2}) is holomorphic 
symmetric bilinear form. Thus it can be written as desired by (\ref{SBF}).
\qed
\vskip2mm

Set
\bea
&&
du={}^t(du_1,...,du_g).
\non
\ena
.

\begin{lemma}\label{lem-oo2}
We have
\bea
&&
\omega(p_1,p_2)=\homega(p_1,p_2)+
{}^tdu(p_1)\eta_1\omega_1^{-1} du(p_2).
\non
\ena
In particular $\eta_1\omega_1^{-1}$ is symmetric.
\end{lemma}
\vskip2mm
\noindent
{\it Proof.}
By Lemma \ref{lem-oo1} 
\bea
&&
\omega-\homega=\sum_{i,j=1}^g {}^tdu(p_1)\,C\,du(p_2),
\label{o-o=cij}
\ena
where $C=(c_{ij})$ is a constant symmetric $g\times g$ matrix.
Since $\int_{\alpha_k}\omega(p_1,p_2)=0$ we have
\bea
&&
\sum_{i=1}^g du_i(p_1)(\eta_1)_{ik}=\sum_{i=1}^g c_{ij}du_i(p_1)(\omega_1)_{jk}.
\non
\ena
Thus
\bea
&&
(\eta_1)_{ik}=\sum_{j=1}^g c_{ij}(\omega_1)_{jk},
\non
\ena
and
\bea
&&
C=\eta_1\omega_1^{-1}.
\non
\ena
\qed
\vskip2mm

\subsection{Symplectic Basis of Cohomology}
For the sake of simplicity we call, hereafter, 
a meromorphic differential on $X$ second kind if it is locally exact. 
In this terminology a first kind differential is of
second kind, the space of differentials of the second kind becomes a vector 
space and the first cohomology group $H^1(X,{\mathbb C})$ is 
described as the space of second kind differentials modulo exact forms.

The intersection
form on $H^1(X,{\mathbb C})$ is given by
\bea
&&
\eta \circ \eta'=\sum \Res \left(\int^p\eta\right)\eta'(p),
\non
\ena
where $\eta$, $\eta'$ are second kind differentials, summation
is over all singular points of $\eta$ and $\eta'$ and $\Res$
means taking a residue at a point.

 Riemann's bilinear relation can be written as
\bea
&&
2\pi i \eta\circ \eta'=
\sum_{i=1}^g \left(
\int_{\alpha_i}\,\eta\int_{\beta_i}\eta'-\int_{\alpha_i}\eta'\int_{\beta_i}\eta
\right).
\non
\ena

\begin{prop}\label{coh-symplectic}
We have
\bea
&&
du_i\circ du_j=0,,
\quad
du_i\circ dr_j=\delta_{ij},
\quad
dr_i\circ dr_j=0,
\label{canonical-coh}
\ena
which means that $\{du_i,dr_j\}$ is a symplectic basis
of $H^1(X,{\mathbb C})$.
\end{prop}
\vskip2mm
\noindent
{\it Proof.}
It is sufficient to prove (\ref{canonical-coh}),
since the linear independence follows from it.

The relation $du_i\circ du_j=0$ is obvious. Let us prove
$du_i\circ dr_j=\delta_{ij}$.   
We calculate $\homega(p_1,p_2)\circ du_j(p_2)$ in two ways.
By Proposition \ref{cijkl} (i)
\bea
\homega(p_1,p_2)\circ du_j(p_2)&=&
\mathop{\Res}_{p_2=p_1}(\int^{p_2}\homega)du_j(p_2)
\non
\\
&=&
-du_j(p_1).
\label{int1}
\ena
On the other hand 
\bea
\homega(p_1,p_2)\circ du_j(p_2)
&=&
\left(
d_{p_2}\Omega(p_1,p_2)+\sum_{i=1}^g du_i(p_1)dr_i(p_2)
\right)
\circ du_j(p_2)
\non
\\
&=&
\sum_{i=1}^g du_i(p_1)(dr_i\circ du_j).
\label{int2}
\ena
Comparing (\ref{int1}) and (\ref{int2}) we have
\bea
&&
dr_i\circ du_j=-\delta_{ij},
\non
\ena
since $\{du_i\}$ are linearly independent.

Next let us prove $dr_i\circ dr_j=0$. Similarly to (\ref{int2}) we have
\bea
&&
\homega(p_1,p_2)\circ dr_j(p_2)
=\sum_{i=1}^g du_i(p_1)(dr_i\circ dr_j).
\label{int3}
\ena
Since $du_k\circ dr_j=\delta_{kj}$ as already proved, we have,
using Lemma \ref{lem-oo2},
\bea
\homega(p_1,p_2)\circ dr_j(p_2)=
\omega(p_1,p_2)\circ dr_j(p_2)
-\sum_{i=1}^g du_i(p_1)(\eta_1\omega_1^{-1})_{ij}.
\label{int4}
\ena
Let us calculate $\omega(p_1,p_2)\circ dr_j(p_2)$.
By Riemann's bilinear relation
\bea
2\pi i \omega(p_1,p_2)\circ dr_j(p_2)
&=&
\sum_{k=1}^g \left(
\int_{\alpha_k}\omega\int_{\beta_k}dr_j-\int_{\alpha_k}dr_j\int_{\beta_k}\omega
\right)
\non
\\
&=&
2\sum_{k=1}^g  (\eta_1)_{jk}\int_{\beta_k}\omega,
\label{int5}
\ena
since $\int_{\alpha_k}\omega=0$.
\vskip2mm

\begin{lemma}\label{int6}
We have
\bea
&&
\int_{\beta_k}\omega=2\pi i\sum_{i=1}^g (2\omega_1)^{-1}_{ki}du_i(p_1),
\label{int7}
\ena
where the integral of the left hand side is with respect to $p_2$
and $(2\omega_1)^{-1}_{ki}$ denotes the $(k,i)$-component of 
$(2\omega_1)^{-1}$.
\end{lemma}
\vskip2mm
\noindent
{\it Proof.}
Similarly to (\ref{int1}) we have
\bea
&&
\omega(p_1,p_2)\circ du_i(p_2)=-du_i(p_1),
\non
\ena
and similarly to (\ref{int5})
\bea
2\pi i\omega(p_1,p_2)\circ du_i(p_2)
&=&
-\sum_{k=1}^g(2\omega_1)_{ik}\int_{\beta_k}\omega.
\non
\ena
The assertion of the lemma follows from these.
\qed
\vskip3mm

Substitute (\ref{int7}) into (\ref{int5}) and get
\bea
&&
\omega(p_1,p_2)\circ dr_j(p_2)
=
\sum_{i=1}^g(\eta_1\omega_1^{-1})_{ji}du_i(p_1).
\label{int8}
\ena
Then we have, by (\ref{int4}),
\bea
&&
\homega(p_1,p_2)\circ dr_j(p_2)=
\sum_{i=1}^g du_i(p_1)
\left(
(\eta_1\omega_1^{-1})_{ji}-(\eta_1\omega_1^{-1})_{ij}
\right),
\non
\ena
which becomes zero since $\eta_1\omega_1^{-1}$ is symmetric by
Lemma \ref{lem-oo2}. It follows from (\ref{int3}) that
$dr_i\circ dr_j=0.$
\qed
\vskip3mm

Due to the relation (\ref{canonical-coh}) Riemann's bilinear
equations take the form
\bea
-{}^t\eta_1\omega_1+{}^t\omega_1\eta_1&=&0,
\label{per-1}
\\
-{}^t\eta_2\omega_2+{}^t\omega_2\eta_2&=&0,
\label{per-2}
\\
-{}^t\eta_1\omega_2+{}^t\omega_1\eta_2&=&-\frac{\pi i}{2}I_g,
\label{per-3}
\ena
where $I_g$ denotes the unit matrix of degree $g$.
If we introduce the matrix
\bea
&&
M=
\left(
\begin{array}{cc}
\omega_1&\omega_2\\
\eta_1&\eta_2\\
\end{array}
\right),
\non
\ena
those relations can be written compactly as
\bea
&&
{}^tM
\left(
\begin{array}{cc}
0&1\\
-1&0\\
\end{array}
\right)
M
=
-\frac{\pi i}{2}
\left(
\begin{array}{cc}
0&I_g\\
-I_g&0\\
\end{array}
\right).
\non
\ena

\section{Schur Function}
Let $p_n(T)$ be the polynomial of $T_1$, $T_2$,... defined by
\bea
&&
\exp(\sum_{n=1}^\infty T_nk^n)=\sum_{n=0}^\infty p_n(T)k^n,
\non
\ena
where $k$ is a variable making a generating function \cite{DJKM}.
\vskip5mm

\noindent
{\bf  Example}
\hskip2mm 
$p_0=1$,
\hskip2mm 
$p_1=T_1$,
\hskip2mm 
$\displaystyle{p_2=T_2+\frac{T_1^2}{2}}$,
\hskip2mm
$\displaystyle{p_3=T_3+T_1T_2+\frac{T_1^3}{6}}$.
\vskip5mm

A sequence of non-negative integers $\lambda=(\lambda_1,...,\lambda_l)$ is 
called a partition if $\lambda_1\geq \cdots\geq \lambda_l$.
We set $|\lambda|=\lambda_1+\cdots+\lambda_l$.
Denote by $\lambda'=(\lambda_1',...,\lambda_{l'}')$, $l'=\lambda_1$,
the conjugate of $\lambda$ \cite{Mac}:
\bea
&&
\lambda_i'=\sharp \{j\,|\lambda_j\geq i\,\}.
\non
\ena

For a partition $\lambda=(\lambda_1,...,\lambda_l)$ define the polynomial
$S_\lambda(T)$ of $T_1$, $T_2$, $T_3$,... by
\bea
&&
S_\lambda(T)=\det (p_{\lambda_i-i+j}(T))_{1\leq i,j\leq l},
\label{s1}
\ena
which we call Schur function.
Notice that, for any $r\geq 0$,
we have
\bea
&&
S_{(\lambda,0^{r})}(T)=S_\lambda(T),
\label{S-extend}
\ena
where $(\lambda,0^{r})=(\lambda_1,...,\lambda_l,0,...,0)$.
\vskip5mm

\noindent
{\bf Example}
\hskip2mm
$S_{(1)}(T)=T_1$, 
\hskip2mm
$\displaystyle{S_{(2,1)}(T)=-T_3+\frac{T_1^3}{3}}$,
\hskip2mm
$\displaystyle{S_{(3,2,1)}(T)=T_1T_5-T_3^2-\frac{1}{3}T_1^3T_3
+\frac{1}{45}T_1^6}$,
\vskip2mm
\hskip17mm
$\displaystyle{S_{(3,1,1)}(T)=T_5-T_1T_2^2+\frac{1}{20}T_1^5}$.
\vskip5mm

We prescribe the degree $-i$ to $T_i$:
\bea
&&
\deg\, T_i=-i.
\non
\ena
The following properties are well known (see for example \cite{DJKM}).
\vskip2mm
\begin{lemma}\label{Schur-prop1}
\vskip2mm
\noindent
(i) $S_\lambda(T)$ is a homogeneous polynomial of degree 
$-|\lambda|$.
\vskip2mm
\noindent
(ii) $S_\lambda(-T)=(-1)^{|\lambda|}S_{\lambda'}(T)$.
\end{lemma}
\vskip2mm

For a partition $\lambda=(\lambda_1,...,\lambda_l)$ we define a
symmetric polynomial of $t_1$, $t_2$, ... $t_l$ by
\bea
&&
s_\lambda(t)=\frac{\det(t_j^{\lambda_i+l-i})_{1\leq i,j\leq l}}
{\prod_{1\leq i<j\leq l}(t_i-t_j)},
\label{s2}
\ena
which we also call Schur function.

Two Schur functions are related by
\bea
&&
S_\lambda(T)=s_\lambda(t) 
\quad
\hbox{if $T_i=\frac{\sum_{j=1}^lt_j^i}{i}$.}
\label{s=S-1}
\ena

If one takes 
$l'\geq |\lambda|$ for $\lambda=(\lambda_1,...,\lambda_l)$,
then the symmetric function $s_{(\lambda,0^{l'-l})}(t)$ 
can be expressed uniquely as a polynomial of power sum symmetric functions 
$T_i=\frac{\sum_{j=1}^{l'}t_j^i}{i}$, $1\leq i\leq l'$. 
This polynomial coincides with $S_\lambda(T)$.

We define a partition associated with a $(n,s)$-curve by
\bea
&&
\lambda(n,s)=(w_g,...,w_1)-(g-1,...,1,0).
\non
\ena
Then

\begin{prop}\label{Schur-prop2}\cite{BEL2}
\vskip2mm
\noindent
(i) $S_{\lambda(n,s)}(T)$ does not depend on the variables other than $T_{w_1}$, ..., $T_{w_g}$, that is, it is a polynomial of the variables $T_{w_1}$, ..., $T_{w_g}$.
\vskip2mm
\noindent
(ii) $\lambda(n,s)'=\lambda(n,s)$.
\vskip2mm
\noindent
(iii) $\displaystyle{|\lambda(n,s)|=\frac{1}{24}(n^2-1)(s^2-1)}$.
\vskip2mm
\noindent
(iv) $\displaystyle{S_{\lambda(n,s)}(-T)=(-1)^{\frac{1}{24}(n^2-1)(s^2-1)}
S_{\lambda(n,s)}(T)}$
\end{prop}
\vskip2mm
Notice that properties (iii) and (iv) follow from (ii) and 
Lemma \ref{Schur-prop1}.

\section{Sigma Function}
\subsection{Definition}
Riemann's constant of a $(n,s)$-curve
with the base point $\infty$ becomes a half period,
since the divisor of the holomorphic one form $du_g$ is $(2g-2)\infty$.
Let
\bea
&&
\qbc{\delta'}{\delta''},
\quad 
\delta',\,\delta''\in \{0,\frac{1}{2}\},
\non
\ena
be the characteristic of Riemann's constant 
$\delta=\delta_0-(g-1)\infty\in J(X)$ 
with respect to our choice $(\infty,\{\alpha_i,\beta_j\})$.
We define the degree of $u_i$ to be $-w_i$:
\bea
&&
\deg\, u_i=-w_i.
\non
\ena
\vskip5mm

\begin{defn}\label{def-1}
The fundamental sigma function or simply the sigma function
$\sigma(u)$ is the holomorphic function on ${\mathbb C}^g$ of the
variables $u={}^t(u_1,...,u_g)$ which satisfies the following
conditions.
\vskip2mm
\noindent
(i) 
\hskip3mm
$
\displaystyle{
\sigma(u+2\omega_1m_1+2\omega_2m_2)/\sigma(u)
=(-1)^{{}^tm_1m_2+2({}^t\delta'm_1-{}^t\delta''m_2)}
}
$
\bea
&&
\qquad\qquad\qquad\qquad
\times
\exp\left(
{}^t(2\eta_1m_1+2\eta_2m_2)(u+\omega_1m_1+\omega_2m_2)
\right).
\label{transformation}
\ena

\vskip2mm
\noindent
(ii) The expansion at the origin takes the form
\bea
&&
\sigma(u)=S_{\lambda(n,s)}(T)|_{T_{w_i}=u_i}+\sum_{d} f_d(u),
\label{normalization}
\ena
where $f_d(u)$ is a homogeneous polynomial of degree $d$ and the range 
of the summation is $d<-|\lambda(n,s)|$.
\end{defn}

It is possible to give an analytic expression of a function satisfying
 the condition (i) in terms of Riemann's theta function.

\begin{prop}\label{sigma-anal}{\rm \cite{BEL1}}
Let $\tau=\omega_1^{-1}\omega_2$.
Then a holomorphic function satisfying (\ref{transformation}) is a constant
multiple of the function
\bea
&&
\exp\left(\frac{1}{2}{}^tu\eta_1\omega_1^{-1}u\right)
\theta\qbc{\delta'}{\delta''}((2\omega_1)^{-1}u,\tau).
\label{analytic-expr}
\ena
\end{prop}

The proposition can easily be proved
using (\ref{per-1}), (\ref{per-2}), (\ref{per-3}) and the uniqueness
of Riemann's theta function $\theta(z)$ \cite{M1}.

The existence of the sigma function is not obvious because of the condition 
(ii). 
In the succeeding subsections we shall construct the sigma function 
explicitly using the algebraic integrals.

\subsection{Algebraic Expression of Prime Form} 
In this section an algebraic expression of the prime form is given
 using the map $x:X\lar {\mathbf P}^1$, $(x,y)\mapsto x$.

Let $B\subset X$ be the set of branch points of the map $x$,
$B'=B\backslash\{\infty\}$,
$\tB=\pi^{-1}(B)\subset \tX$, $\tB'=\pi^{-1}(B')$
$\bB=x(B)\subset {\mathbb P}^1$ and $\bB'=x(B')$.
Let $\tilde{p}\in\tX$ and $\tgamma$ be a path from $\tinfty$ to $\tilde{p}$
in $\tX\backslash \tB'$. Then $\bgamma=(x\circ \pi)(\tgamma)$ is a path
from $\infty$ to $(x\circ \pi)(\tilde{p})$. Let $\igamma$ be the lift of 
$\bgamma$ to $X\backslash B'$ connecting $p^{(i)}$ to $\infty$ and
$\tigamma$ the lift of $\igamma$ to $\tX$ beginning at $\tinfty$.
Denote $\tilde{p}^{(i)}$ the end point of $\tigamma$. 
Then $\tilde{p}^{(i)}$ lies over $\tip$.

In this way for each path $\tgamma=\tgamma^{(0)}$ from $\tinfty$
to $\tilde{p}$ in $\tX\backslash \tB'$ we have a uniquely determined path
from $\tinfty$ to $\tilde{p}^{(i)}$ in $\tX\backslash \tB'$.
Let $\tgamma_i=\tgamma_i^{(0)}$, $i=1,2$ be paths from $\tinfty$ to $\tilde{p}_i$ and
$\tgamma_i^{(j)}$ be the corresponding path from $\tinfty$ to $\tilde{p}_i^{(j)}$.
The path from $\tilde{p}_1^{(i)}$ to $\tilde{p}_2^{(i)}$ is defined by
$\tgamma_2^{(i)}\circ (\tgamma_1^{(i)})^{-1}$.
Hereafter, for $\tilde{p}\in X$, we denote $\pi(\tilde{p})$ by $p$ 
if there is no fear of confusion.

\begin{prop}\label{alg-pf}\cite{F}
We have
\bea
&&
E(\tilde{p}_1,\tilde{p}_2)^2=\frac{(x(p_2)-x(p_1))^2}{dx(p_1)dx(p_2)}
\exp\left(
\sum_{i=1}^{n-1}\int_{\tilde{p}_1^{(i)}}^{\tilde{p}_2^{(i)}}\int_{\tilde{p}_1}^{\tilde{p}_2}
\omega
\right).
\non
\ena
\end{prop}
\vskip8mm

\noindent
{\bf Remark} 
In Fay's book \cite{F} (p17, formula (v)), the prime form is expressed 
by the right hand side of Proposition \ref{alg-pf} multiplied by
the term $\exp(\sum_{i=1}^g\int_{p_1}^{p_2}m_idv_i)$.
This difference stems from the fact that $E(\tilde{p}_1,\tilde{p}_2)$ is described
on the universal covering $\tX$ in this paper while it is described
in the fundamental polygon cut out along the homology basis 
$\{\alpha_i,\beta_i\}$ in \cite{F}. 
\vskip5mm

For the sake to be complete and self-contained
we give a proof of this proposition.

\begin{lemma}\label{lem-rational}
We have
\bea
&&
\frac{(x(w)-x(p_2))(x(z)-x(p_1))}{(x(w)-x(p_1))(x(z)-x(p_2))}
=\exp\left(
\sum_{i=0}^{n-1}\int_{\tilde{z}}^{\tilde{w}}
\int_{\tilde{p}_1^{(i)}}^{\tilde{p}_2^{(i)}}
\omega
\right).
\label{f-rational}
\ena
\end{lemma}
\vskip2mm

\noindent
{\it Proof.}
By Proposition \ref{exp-omega} 
\bea
&&
\exp\left(
\sum_{i=0}^{n-1}\int_{\tilde{z}}^{\tilde{w}}
\int_{\tilde{p}_1^{(i)}}^{\tilde{p}_2^{(i)}}
\omega
\right)
=\prod_{i=0}^{n-1}
\frac{E(\tilde{w},\tilde{p}_2^{(i)})E(\tilde{z},\tilde{p}_1^{(i)})}
{E(\tilde{w},\tilde{p}_1^{(i)})E(\tilde{z},\tilde{p}_2^{(i)})}.
\non
\ena
Let us consider the right hand side of this equation as a function 
of $\tilde{w}$ and
denote it by $F(\tilde{w})$. 
By the property (iv) of the prime form,
if the abelian image of $\gamma\in \pi_1(X,\infty)$ 
is $\sum_{i=1}^g m_{1,i}\alpha_i+\sum_{i=1}^g m_{2,i}\beta_i$ then
\bea
&&
F(\gamma \tilde{w})=F(\tilde{w})
\exp\left(
-2\pi i\sum_{j=0}^{n-1}{}^tm_2\int_{\tilde{p}_1^{(j)}}^{\tilde{p}_2^{(j)}}dv
\right).
\label{F-transf}
\ena
\vskip5mm

\begin{lemma}\label{lem-integral}
\bea
&&
\sum_{j=0}^{n-1}\int_{\tilde{p}_1^{(j)}}^{\tilde{p}_2^{(j)}}dv_i=0,
\quad
1\leq i\leq g.
\non
\ena
\end{lemma}
\vskip2mm
\noindent
{\it Proof.}
It is sufficient to prove
\bea
&&
\sum_{j=0}^{n-1}\int_{\tilde{p}_1^{(j)}}^{\tilde{p}_2^{(j)}}du_i=0.
\label{int-dui}
\ena
Let us fix $\tilde{p}_1$ and consider the left hand side of (\ref{int-dui})
as a function of $\tilde{p}_2$. We denote it by $G(\tilde{p}_2)$.
Then $G(\tilde{p}_2)$ is a holomorphic function on 
$Y=\tX-\tB$. At each $\tilde{p}^{(i)}\in Y$
one can take $x$-coordinate as a local coordinate around it.
By differentiating $G(\tilde{p_2})$ with respect to the local coordinate $x$
we get
\bea
&&
dG(\tilde{p_2})=\sum_{j=0}^{n-1} \frac{x^{a_i-1}(y^{(j)})^{n-1-b_i}}{f_y(x,y^{(j)})}
dx,
\label{dg}
\ena
where $p_2^{(j)}=(x,y^{(j)})$. 
\vskip5mm

\begin{lemma}\label{res-sum}
Let $q\in X$ and $q^{(i)}=(x,y^{(i)})$.
Suppose that $q\notin B$.
Then
\bea
&&
\sum_{j=0}^{n-1}\frac{(y^{(j)})^b}{f_y(x,y^{(j)})}=0,
\qquad
0\leq b\leq n-2.
\non
\ena
\end{lemma}
\vskip2mm
\noindent
{\it Proof.}
Let 
\bea
&&
f(x,y)=\prod_{j=0}^{n-1}(y-y^{(j)}).
\non
\ena
Since $q$ is not a branch point, $y^{(i)}\neq y^{(j)}$, $i\neq j$.
For $0\leq b\leq n-2$ we have
\bea
&&
\mathop{\Res}_{z=\infty}\,\,\frac{z^b}{\prod_{j=0}^{n-1}(z-y^{(j)})}\,dz=0,
\non
\ena
where $z$ is a variable on ${\mathbb P}^1$.
Thus, by the residue theorem on ${\mathbb P}^1$, 
\bea
&&
\sum_{j=0}^{n-1}
\mathop{\Res}_{z=y^{(j)}}\,\,\frac{z^b}{\prod_{i=0}^{n-1}(z-y^{(i)})}\,dz
=\sum_{j=0}^{n-1}
\frac{(y^{(j)})^b}{f_y(x,y^{(j)})}
=0.
\non
\ena
\qed
\vskip5mm

By Lemma \ref{res-sum}  and (\ref{dg}) we have
\bea
&&
dG(\tilde{p_2})=0,
\non
\ena
on $Y$. Since $G(\tilde{p_2})$ is continuous at each point of $\tB$,
it is a constant on $\tX$. By the definition of the integration path
$G(\tilde{p_1})=0$.
Therefore $G(\tilde{p_2})$ is identically zero as desired. \qed
\vskip3mm

Let us continue the proof of Lemma \ref{lem-rational}.
By Lemma \ref{lem-integral} and (\ref{F-transf}) $F(\tilde{w})$ is 
$\pi_1(X,\infty)$-invariant and can be considered as a meromorphic 
function on $X$.
By comparing zeros and poles,
\bea
&&
F(\tilde{w})=C\frac{x(w)-x(p_2)}{x(w)-x(p_1)},
\non
\ena
for some constant $C$. Since $F(\tilde{z})=1$ 
\bea
&&
C=\frac{x(z)-x(p_1)}{x(z)-x(p_2)}
\non
\ena
which proves the lemma.
\qed
\vskip5mm

\noindent
Proof of Proposition \ref{alg-pf}.
\vskip2mm

\noindent
In Lemma \ref{lem-rational} take the limit 
$\tilde{z}\rightarrow \tilde{p}_1$, 
$\tilde{w}\rightarrow \tilde{p}_2$ and use
\bea
&&
\lim_{\tilde{w}\rightarrow \tilde{q}}\frac{x(w)-x(q)}{E(\tilde{w},\tilde{q})}=-dx(q),
\non
\\
&&
\exp\left(
\int_{\tilde{z}}^{\tilde{w}}\int_{\tilde{p}_1}^{\tilde{p}_2}
\omega
\right)
=
\frac{E(\tilde{w},\tilde{p}_2)E(\tilde{z},\tilde{p}_1)}
{E(\tilde{w},\tilde{p}_1)E(\tilde{z},\tilde{p}_2)}.
\ena
Then we easily get the desired result.
\qed

\subsection{Prime Function}
Let $\sqrt{du_g}$ be the holomorphic section of the line bundle on $X$
defined by the divisor $(g-1)\infty$ satisfying
\bea
(\sqrt{du_g})^2&=&du_g,
\label{normal-1}
\\
\sqrt{du_g}&=&t^{g-1}\left(1+O(t)\right)\sqrt{dt},
\label{normal-2}
\ena
where $t$ is the local parameter (\ref{par-inf}) around $\infty$.

\vskip5mm
\begin{defn}\label{prime-func}
We define the prime function $\tE(\tilde{p}_1,\tilde{p}_2)$ on
$\tX\times\tX$ by
\bea
&&
\!\!\!\!\!\!\!\!\!\!\!\!
\tE(\tilde{p}_1,\tilde{p}_2)=-E(\tilde{p}_1,\tilde{p}_2)
\sqrt{du_g(p_1)}
\sqrt{du_g(p_2)}
\exp\left(\frac{1}{2}
\int_{\tilde{p}_1}^{\tilde{p}_2}{}^t du \,\cdot\,\eta_1\omega_1^{-1} \cdot\,
\int_{\tilde{p}_1}^{\tilde{p}_2} du
\right).
\label{def-prf}
\ena
\vskip2mm
\end{defn}
\vskip3mm

Since 
\bea
&&
\delta=\delta'\tau+\delta''=\delta_0-(g-1)\infty=(g-1)\infty-\delta_0
\quad
\hbox{in}
\quad
J(X),
\non
\ena
$\tE(\tilde{p}_1,\tilde{p}_2)$ can be considered as a holomorphic section of the line bundle 
$
\pi_1^\ast {\cal L}_{\delta}\otimes 
\pi_2^\ast {\cal L}_\delta\otimes I_2^\ast\Theta
$
on $X\times X$.

By Proposition \ref{alg-pf} and Lemma \ref{lem-oo2} we have
\vskip2mm
\bea
&&
\tE(\tilde{p}_1,\tilde{p}_2)^2=
\frac{\left(x(p_2)-x(p_1)\right)^2}{f_y(p_1)f_y(p_2)}
\exp\left(
\sum_{i=1}^{n-1}\int_{\tilde{p}_1^{(i)}}^{\tilde{p}_2^{(i)}}
\int_{\tilde{p}_1}^{\tilde{p}_2}
\homega
\right).
\label{alg-pfn}
\ena
\vskip2mm

We need to put one of the variables in 
$\tE(\tilde{p}_1,\tilde{p}_2)$
to be $\tinfty$ in order to describe the sigma function.
Since $\tE(\tilde{p}_1,\tilde{p}_2)$ becomes zero at $\tilde{p}_i=\tinfty$,
it is defined in the following manner.
Take the local coordinate $t$ (\ref{par-inf}) and the local frame $\sqrt{dt}$
as above and define
\bea
&&
\!\!\!\!\!\!\!\!\!\!\!\!\!\!\!\!\!\!\!\!\!\!\!\!\!\!\!\!\!\!
\!\!\!\!\!\!\!\!\!\!\!\!\!\!\!\!\!\!\!\!\!\!\!\!\!
E(\tinfty,\tilde{p}_2)=
E(\tilde{p}_1,\tilde{p}_2)\sqrt{dt(p_1)}\vert_{t(p_1)=0},
\label{inf-pf}
\ena
\bea
&&
\qquad\quad
\tE(\tinfty,\tilde{p})=
E(\tinfty,\tilde{p})
\sqrt{du_g(p)}
\exp\left(\frac{1}{2}
\int_{\tinfty}^{\tilde{p}}{}^t du \,\cdot\,\eta_1\omega_1^{-1} \cdot\,
\int_{\tinfty}^{\tilde{p}} du
\right).
\label{inf-pfn}
\ena
\vskip2mm

Notice that $E(\tinfty,\tilde{p})$ and $\tE(\tinfty,\tilde{p}_2)$
can be considered as holomorphic sections of 
$L_0^{-1}\otimes I_1^\ast \Theta$ and 
${\cal L}_\delta\otimes I_1^\ast \Theta$
respectively.
By (\ref{normal-2}) and the property (iii) of $E(\tilde{p}_1,\tilde{p}_2)$
we have
\bea
&&
-\tE(\tilde{p}_1,\tilde{p}_2)\sqrt{dt(p_1)}=
\tE(\tinfty,\tilde{p}_2)t(p_1)^{g-1}
+
O\left(t(p_1)^g\right).
\non
\ena
\vskip2mm
The properties of the prime form imply those
of $\tE(\tilde{p}_1,\tilde{p}_2)$ and $\tE(\tinfty,\tilde{p})$.
The next proposition follows from the properties (i) and 
(ii) of the prime form.
\vskip3mm

\begin{prop}\label{prop-tE}
(i) $\tE(\tilde{p}_2,\tilde{p}_1)=-\tE(\tilde{p}_1,\tilde{p}_2)$.
\vskip2mm
\noindent
(ii) As a section of a line bundle on $X\times X$, the zero divisor
of $\tE(\tilde{p}_1,\tilde{p}_2)$  is 
\bea
&&
\Delta+(g-1)(\{\infty\}\times X+X\times \{\infty\}).
\non
\ena
\vskip2mm
\noindent
(iii) As a section of a line bundle on $X$ the zero divisor of 
$\tE(\tinfty,\tilde{p})$  is $g\infty$.
\end{prop}

Later we shall study the series expansion of those functions 
(see Lemma \ref{te-exp}).

\begin{prop}\label{te-transf}
Let the abelian image of $\gamma\in \pi_1(X,\infty)$ be
 $\sum m_{1,i}\alpha_i+\sum m_{2,i}\beta_i$.
Then
\vskip5mm
\noindent
(i)
\hskip3mm
$
\displaystyle{
\tE(\tilde{p}_1,\gamma \tilde{p}_2)/\tE(\tilde{p}_1,\tilde{p}_2)
=(-1)^{{}^tm_1m_2+2({}^t\delta'm_1-{}^t\delta''m_2)}
}
$
\bea
&&
\qquad\quad
\times
\exp\left(
{}^t(2\eta_1m_1+2\eta_2m_2)(\int_{\tilde{p}_1}^{\tilde{p}_2}du+\omega_1m_1+\omega_2m_2)
\right).
\non
\ena
\vskip2mm
\noindent
(ii) The equation (i) substituted by $\tilde{p}_1=\tinfty$ holds for 
$\tE(\tinfty,\tilde{p}_2)$.
\end{prop}
\vskip2mm
\noindent
{\it Proof.} (i)
Consider
\bea
&&
F_1(\tilde{p}_1,\tilde{p}_2)=
E(\tilde{p}_1,\tilde{p}_2)
\sqrt{du_g(p_1)du_g(p_2)}.
\non
\ena
It is a section of the bundle $\pi_1^\ast {\cal L}_{\delta}\otimes 
\pi_2^\ast {\cal L}_{\delta}\otimes \pi_{12}^\ast \Theta$. 
For a non-singular odd half period $\alpha=\tau \alpha'+\alpha''$
set
$$
F_2=F_1/\theta[\alpha](\int_{\tilde{p}_1}^{\tilde{p}_2} dv),
$$
which is a section of the line bundle 
$\pi_1^\ast {\cal L}_{\delta-\alpha}\otimes 
\pi_2^\ast {\cal L}_{\delta-\alpha}$.
Let 
\bea
&&
F_2(\tilde{p}_1,\gamma\tilde{p}_2)=\chi(\gamma)F_2(\tilde{p}_1,\tilde{p}_2),
\quad
\gamma\in \pi_1(X,\infty),
\non
\ena
where $\chi:\pi_1(X,\infty)\lar {\mathbb C}^\ast$ is a representation
of $\pi_1(\infty,X)$.
Since
\bea
&&
\frac{du_g}{h_\alpha^2}
=\frac{1}{\sum_{i,j=1}^g \frac{\partial \theta[\alpha]}{\partial z_i}(0)
(2\omega_1)^{-1}_{ij}x^{a_j-1}y^{n-1-b_j}}
\non
\ena
is $\pi_1(X,\infty)$-invariant, so is $F_2^2$. Thus $\chi(\gamma)^2=1$ and
$\chi$ is a unitary representation. Therefore, if the abelian image of 
$\gamma$ is $\sum m_{1,i}\alpha_i+\sum m_{2,i}\beta_i$, we have
\bea
&&
\chi(\gamma)=
\exp\left(
2\pi i({}^t (\delta'-\alpha')m_1+{}^t(\delta''-\alpha'')m_2)
\right),
\non
\ena
by (\ref{unitary}) and consequently
\bea
&&
\frac{F_1(\tilde{p}_1,\gamma \tilde{p}_2)}{F_1(\tilde{p}_1,\tilde{p}_2)}
=(-1)^{2({}^t\delta'm_1-{}^t\delta''m_2)}
\exp\left(
-\pi i {}^tm_2\tau m_2-2\pi i{}^t m_2\int_{\tilde{p}_1}^{\tilde{p}_2}dv
\right),
\label{F1-transf}
\ena
by (\ref{R-theta1}).

Next let 
\bea
&&
F_3(\tilde{p}_1,\tilde{p}_2)=
\exp\left(\frac{1}{2}
\int_{\tilde{p}_1}^{\tilde{p}_2}du\,\,\, \eta_1\omega_1^{-1}\! \int_{\tilde{p}_1}^{\tilde{p}_2}du
\right).
\ena
By calculations we have
\vskip5mm
\hskip5mm
$
\displaystyle{
\frac{F_3(\tilde{p}_1,\gamma \tilde{p}_2)}{F_3(\tilde{p}_1,\tilde{p}_2)}
=
\exp\left(
{}^t(2\eta_1m_1+2\eta_2m_2)(\int_{\tilde{p}_1}^{\tilde{p}_2}du+\omega_1m_1+\omega_2m_2)
+\pi i {}^tm_1m_2\right)
}
$
\bea
&&
\!\!\!\!\!\!\!\!\!\!\!\!\!\!\!\!\!\!\!\!\!\!\!\!\!\!\!\!\!\!
\times
\exp\left(
\pi i {}^t m_2\tau m_2+2\pi i{}^tm_2\int_{\tilde{p}_1}^{\tilde{p}_2}dv
\right).
\label{F2-transf}
\ena
Here we use (\ref{per-1}), (\ref{per-3}) and the relation 
\bea
&&
dv=(2\omega_1)^{-1}du.
\non
\ena
Multiplying (\ref{F1-transf}) and (\ref{F2-transf}) we get the desired result.
\vskip2mm

\noindent
(ii) The statement is obvious from (\ref{inf-pf}), (\ref{inf-pfn}) and (i).
\qed

\subsection{Algebraic Expression of Sigma Function}
We now state our main theorems of this paper.
\begin{theorem}\label{sigma-1}
For $N\geq g$
\bea
&&
\sigma(\,\,\sum_{i=1}^N\int_{\tinfty}^{\tilde{p}_i}du\,\,)\,\,
=
\,\,
\frac{\prod_{i=1}^N\tE(\tinfty,\tilde{p}_i)^N}{\prod_{i<j}\tE(\tilde{p}_i,\tilde{p}_j)}
\,\,\det(f_i(p_j))_{1\leq i,j\leq N}.
\label{main-1}
\ena
\end{theorem}
\vskip5mm
\noindent

It is possible to derive a more general formula which 
contains Theorem \ref{sigma-1} as a limit.

Let $N\geq g$, $p_i$, $q_i$, $i=1,...,N$ be points on $X$ and
$f_i$, $i=1,...,nN$ be the basis of $L((nN+g-1)\infty)$ defined before.
Consider the function
\vskip8mm
\hskip5mm
$
\displaystyle{
F_N=
\frac{D_N}{\prod_{i<j}(x(q_i)-x(q_j))^{n-2}
\prod_{k=1}^{N}\prod_{1\leq i<j\leq n-1}
\left(y(q_k^{(i)})-y(q_k^{(j)})\right)},
}
$
\vskip10mm
\hskip5mm
$
\displaystyle{
D_N=
\left|
\begin{array}{cccccccccc}
\!\!f_1(p_1)\!\!&\!\!\cdots\!\!&\!\!f_1(p_N)\!\!&\!\!f_1(q_1^{(1)})\!\!&
\!\!\cdots\!\!&\!\!f_1(q_1^{(n-1)})\!\!&
\!\!\cdots\!\!&\!\!f_1(q_N^{(1)})\!\!&\!\!\cdots\!\!&
\!\!f_1(q_N^{(n-1)})\!\!
\\
\!\!\vdots&\!\!\quad\!\!&\!\!\vdots\!\!&\!\!\vdots\!\!&\!\!\quad\!\!&
\!\!\vdots\!\!&
\!\!\quad\!\!&\!\!\vdots\!\!&\!\!\quad\!\!&\!\!\vdots\!\!\\
\!\!f_{nN}(p_1)\!&\!\!\cdots\!\!&\!f_{nN}(p_N)\!&\!f_{nN}(q_1^{(1)})\!&
\!\!\cdots\!\!&\!\!f_{nN}(q_1^{(n-1)})\!\!&
\!\!\cdots\!\!&\!\!f_{nN}(q_N^{(1)})\!\!&\!\!\cdots\!\!&
\!\!f_{nN}(q_N^{(n-1)})\!\!
\\
\end{array}
\right|.
}
$
\vskip7mm

Notice that, for each $j$, $F_N$ is symmetric in $q_j^{(1)},...,q_j^{(n-1)}$
and does not depend on the way of labeling the points 
$x^{-1}(x(q_j))\backslash \{q_j\}$.
\vskip5mm

\begin{theorem}\label{sigma-2}
Suppose that $N\geq g$. Then
\bea
&&
\sigma(\,\,\sum_{i=1}^N\int_{\tilq_{i}}^{\tilde{p}_i}du\,\,)\,\,
=\,\,C_N\,\,
M_N
F_N,
\label{general}
\ena
where
\vskip5mm
\hskip12mm
$
\displaystyle{
M_N=\frac{\prod_{i,j=1}^N\tE(\tilde{p}_i,\tilq_j)}
{\prod_{i<j}\left(\tE(\tilde{p}_i,\tilde{p}_j)\tE(\tilq_i,\tilq_j)\right)
\prod_{i,j=1}^N\left(x(p_i)-x(q_j)\right)},
}
$
\vskip5mm
\bea
&&
C_N=(-1)^{\frac{1}{2}nN(N-1)}
\left(\frac{\epsilon(s)}{\epsilon(1)}\right)^N
\epsilon_n^{\frac{1}{2}N(N-1)-\frac{1}{4}N(N-1)(n-1)(n-2)+\frac{1}{2}Nn(n-1)
-\frac{1}{2}gNn(n-3)},
\label{cn}
\ena
\vskip5mm
\hskip15mm
$
\epsilon_n=\exp(2\pi i/n),
\qquad
\epsilon(r)=\prod_{1\leq i<j\leq n-1}(\epsilon_n^{ri}-\epsilon_n^{rj}).
$
\end{theorem}
\vskip5mm

\noindent
{\it Proof of Theorem \ref{sigma-1}}
\vskip2mm
Let $G$ be the right hand side of (\ref{main-1}) divided by
the left hand side.
Obviously $G$ is a symmetric function of $\tilde{p}_1,...,\tilde{p}_N$.
Using Proposition \ref{te-transf} and Proposition \ref{prop-tE} one can easily
verify the following properties.
\vskip3mm
\noindent
(i) $G(\gamma \tilde{p}_1,\tilde{p}_2,...,\tilde{p}_N)=G(\tilde{p}_1,\tilde{p}_2,...,\tilde{p}_N)$
for any $\gamma\in \pi_1(X,\infty)$.
\vskip2mm
\noindent
(ii) The right hand side of (\ref{main-1}) is holomorphic as a function of $\tilde{p}_1$.
\vskip3mm

Let us consider $G$ as a function of $\tilde{p}_1,...,\tilde{p}_g$.
By (i) it can be considered as a meromorphic function on the
$g$-th symmetric product $S^gX=X^g/S_g$ and therefore on the Jacobian 
$J(X)$. As a meromorphic function on $J(X)$ $G$ has poles only on 
$\Sigma=\{\sigma(u)=0\}$ of order at most one by (ii). Thus it is a constant
which means that it is independent of $\tilde{p}_i$, $1\leq i\leq g$.
Since $G$ is symmetric, it does not depend on all of the variables.
The constant is calculated in the proof of Theorem \ref{main-cor} (i) in
the next section.
\qed
\vskip5mm

In order to prove Theorem \ref{sigma-2} we have to study the properties
of the function $F_N$.

\begin{lemma}\label{lem-generalized}
(i) $F_N$ is skew symmetric with respect to $\{p_i\}$ and $\{q_i\}$ respectively.
\vskip2mm
\noindent
(ii) In each of the variables $\{p_i,q_j\}$ $F_N$ is a meromorphic function .
\vskip2mm
\noindent
(iii) As a function of $p_1$  $F_N$ has poles only at $\infty$ of order
 at most $nN+g-1$ and  zeros at 
$p_j$ $(j\geq 2)$, $q_k^{(i)}$ $(1\leq k\leq N, 1\leq i\leq n-1)$.
\vskip2mm
\noindent
(iv) As a function of $q_1$  $F_N$ has poles only at $\infty$ of order
 at most $nN+g-1$ and  zeros at 
$q_j$ $(j\geq 2)$, $p_k^{(i)}$ $(1\leq k\leq N, 1\leq i\leq n-1)$.
\end{lemma}
\vskip2mm
{\it Proof.}
(i) The skew symmetry in $\{p_i\}$ is obvious and that in $\{q_i\}$ 
can be easily verified.
\vskip2mm
\noindent
(ii) It is obvious that $F_N$ is a meromorphic function with respect
to $p_i$. Let us prove that $F_N$ is a meromorphic function of $q_1$.
To this end we first prove that $F_N$ is a symmetric polynomial
of $y(q_1^{(1)})$,...,$y(q_1^{(n-1)})$. 
Notice that $D_N$ is
a skew symmetric polynomial of $y(q_1^{(1)})$,...,$y(q_1^{(n-1)})$,
since $f_k(q_1^{(i)})=x(q_1)^ay(q_1^{(i)})^b$ for some $a,b$ and
permuting $q_1^{(1)},...,q_1^{(n-1)}$ is the same as permuting columns.
Thus 
\bea
&&
D_N/\prod_{i<j}\left(y(q_1^{(i)})-y(q_1^{(j)})\right),
\label{divide-DN}
\ena
is a symmetric polynomial of $y(q_1^{(1)})$,...,$y(q_1^{(n-1)})$.
Obviously its coefficients are polynomials of $x(q_1)$.

Now it is sufficient to prove that any symmetric polynomial
of $y(q^{(1)})$,...,$y(q^{(n-1)})$ is a polynomial of $x(q)$ and $y(q)$.
Let us write
\bea
&&
f(x(q),y)=\sum_{i=0}^n A_i(x(q)) y^{n-i}=\prod_{i=0}^\infty (y-y_i),
\non
\ena
where $A_0=1$, $A_i(x)$ is a polynomial of $x$ and $y_i=y(q^{(i)})$.
Then 
\bea
&&
e_i(y_0,...,y_{n-1})=(-1)^iA_i(x(q)),
\non
\ena
where $e_i(t_1,...,t_n)$ is the $i$-th elementary symmetric polynomial,
\bea
&&
\prod_{i=1}^n(t+t_i)=\sum_{i=0}^n e_i(t_1,...,t_n)t^{n-i}.
\non
\ena
Using
\bea
&&
e_i(y_0,...,y_{n-1})=y_0e_{i-1}(y_1,...,y_{n-1})+e_i(y_1,...,y_{n-1}),
\non
\ena
one can easily prove that every $e_i(y_1,...,y_{n-1})$ is a polynomial
of $y_0=y(q)$ and $x(q)$. Thus any symmetric polynomial of $y_1,...,y_{n-1}$
is a polynomial of $x(q)$ and $y(q)$.
\vskip2mm
\noindent
(iii) This is obvious.
\vskip2mm
\noindent
(iv) As proved in (ii) 
(\ref{divide-DN}) is a polynomial of $x(q_1)$, $y(q_1)$ and therefore its
only singularity is $\infty$. 
Let us examine the zeros of $D_N$.

Notice that $D_N$ has zeros at $q_1=q_j$, $j\neq 1$ of
order at least $n-1$. In fact $q_1=q_j$ implies that
$\{q_1^{(1)},...,q_1^{(n-1)}\}=\{q_j^{(1)},...,q_j^{(n-1)}\}$.
Therefore $D_N$ has zeros of order at least $n-1$ at $q_1=q_j$.

In a similar manner $D_N$ has zeros at $q_1=q_j^{(i)}$, $i\neq 0$,
of order at least 
$n-2$, because, in this case, the number of elements in 
$\{q_1^{(1)},...,q_1^{(n-1)}\}\cap \{q_j^{(1)},...,q_j^{(n-1)}\}$
is $n-2$.

Consequently the only singularity of $F_N$ is at most $q_1=\infty$.
Moreover $F_N$ has zeros at $q_1=q_j$, $j\neq 1$. The fact that
$F_N$ is zero at $q_1=p_j^{(i)}$, $i\neq 0$ can be easily verified.

The order of poles at $\infty$, which we denote by $\ord_\infty\, F_N$,
is estimated as
\bea
\ord_\infty\,F_N&\leq& O_1-O_2-O_3,
\non
\\
O_1&=&(nN+g-1)+\cdots+\left(nN+g-(n-2)\right),
\non
\\
O_2&=&(n-2)n(N-1),
\non
\\
O_3&=&s\bc{n-1}{2},
\non
\ena
where $O_1$ is the maximal possible order of poles of $D_N$,
$O_2$ is the order of poles of $\prod_{j=2}^N(x(q_1)-x(q_j))^{n-2}$
and $O_3$ is that of $\prod_{i<j}(y(q_1^{(i)})-y(q_1^{(j)}))$.
By calculation we get
\bea
&&
O_1-O_2-O_3=nN+g-1.
\non
\ena
\qed

\vskip5mm
\noindent
{\it Proof of Theorem \ref{sigma-2}.}
The proof is similar to that of Theorem \ref{sigma-1}.
Let $G(p_1,...,p_N|q_1,...,q_N)$ be the right hand side divided
by the left hand side of (\ref{general}).

The function $G$ is $\pi_1(X,\infty)$-invariant in each 
of the variables $\{\tilde{p}_i,\tilq_j\}$ by Proposition \ref{te-transf}.
Consider $G$ as a function of $p_1,...,p_g$. 
By Lemma \ref{lem-generalized} $G$ can be considered 
as a meromorphic function on $J(X)$ which has poles only on $\Sigma$
of order at most one.
Thus it does not depend on $p_1,...,p_g$. Since $G$ is symmetric
in $\{p_i\}$ it does not depend on any $p_i$. Similarly $G$ does not
depend on $\{q_i\}$. Thus $G$ is a constant. The constant is calculated
in the next section.
\qed
\vskip5mm

\subsection{The Series Expansion of Sigma Function}
In this section we study the series expansion of the sigma function
using the expression obtained in the previous section.

Let us use the multi-index notations like
\bea
&&
u^{\alpha}=u_{1}^{\alpha_1}\cdots u_{g}^{\alpha_g},
\qquad
\alpha=(\alpha_1,...,\alpha_g).
\non
\ena
We set
\bea
&&
|\alpha|=\sum_{i=1}^gw_i \alpha_i,
\non
\ena
so that $\deg\, u^\alpha=-|\alpha|$.

\begin{theorem}\label{main-cor}
(i) The expansion of $\sigma(u)$
 at the origin takes the form
\bea
&&
\sigma(u)=S_{\lambda(n,s)}(T)|_{T_{w_i}=u_i}+\sum a_\alpha u^\alpha,
\non
\ena
where $a_\alpha\in {\mathbb Q}[\{\lambda_{ij}\}]$ and the sum
is taken for $\alpha$ such that $|\alpha|>|\lambda(n,s)|$.
\vskip2mm
\noindent
(ii) In (i) $a_\alpha$ is homogeneous of degree 
$-|\lambda(n,s)|+|\alpha|$. In particular $\sigma(u)$ is homogeneous of degree
$-|\lambda(n,s)|$ with respect to the variables $\{u_i,\lambda_{jk}\}$.
\vskip2mm
\noindent
(iii) $\sigma(-u)=(-1)^{|\lambda(n,s)|}\sigma(u)$.
\end{theorem}

In the remaining of this section $t$ denotes the local parameter 
around $\infty$ specified by (\ref{par-inf}).
\bea
&&
\deg\, t=-1.
\non
\ena
\begin{lemma}\label{lem-expansion}
(i) If we write $y=\frac{1}{t^s}\sum_{i=0}^\infty c_it^i$, $c_0=1$,
then $c_i$ belongs to ${\mathbb Q}[\{\lambda_{kl}\}]$ and is homogeneous
of degree $i$.
\vskip2mm
\noindent
(ii) If we write $dx/f_y=-t^{2g-2}(1+\sum_{i=1}^\infty c'_it^i)dt$, then 
$c'_i$ belongs to ${\mathbb Q}[\{\lambda_{kl}\}]$ and is homogeneous
of degree $i$.
\vskip2mm
\noindent
(iii) Let $t_i=t(p_i)$ and
\bea
&&
\homega(p_1,p_2)=
\left(
\frac{1}{(t_1-t_2)^2}
+\sum_{i,j=0}^\infty a_{ij}t_1^it_2^j
\right)dt_1dt_2.
\non
\ena
Then $a_{ij}=a_{ji}$, $a_{ij}$ belongs to 
${\mathbb Q}[\{\lambda_{kl}\}]$ and is homogeneous
of degree $i+j+2$.
\end{lemma}
\vskip2mm
\noindent
{\it Proof.} (i)
Substitute the expressions of $x$, $y$ in $t$ to the defining equation
$f(x,y)=0$ with $f$ being given by (\ref{ns-curve}). Then 
\bea
&&
y=\frac{1}{t^s}\left(
1+\sum_{ni+sj<ns}\lambda_{ij}t^{ns-ni-sj}(\sum_{r=0}^\infty c_rt^r)^j
\right)^{\frac{1}{n}}.
\non
\ena
If we write
\bea
&&
1+\sum_{ni+sj<ns}\lambda_{ij}t^{ns-ni-sj}(\sum_{r=0}^\infty c_rt^r)^j
=\sum_{i=0}^\infty b_i t^i,
\quad
b_0=1,
\label{y-1}
\ena
then each $b_i$ is a polynomial of $c_1,...,c_{i-1}$ and $\{\lambda_{kl}\}$
with the coefficients in ${\mathbb Q}$, since $ns-ni-sj>0$.

In the expansion
\bea
&&
(\sum_{i=0}^\infty b_it^i)^{1/n}=
1+\sum_{j=1}^\infty \bc{1/n}{j}(\sum_{i=1}^\infty b_it^i)^j,
\label{y-2}
\ena
the coefficient of $t^r$ is a polynomial of $b_1,...,b_r$ 
with the coefficient in ${\mathbb Q}$ and consequently 
is a polynomial of $c_1,...,c_{r-1}$, 
$\{\lambda_{kl}\}$ with the coefficient in ${\mathbb Q}$.
Thus the equation
\bea
&&
\sum_{i=0}^\infty c_it^i=(\sum_{i=0}^\infty b_i t^i)^{1/n}
\label{y-3}
\ena
implies that $c_r\in {\mathbb Q}[c_1,...,c_{r-1},\{\lambda_{kl}\}]$ for
$r\geq 1$. This proves $c_r\in {\mathbb Q}[\{\lambda_{kl}\}]$.

Next let us prove $\deg\, c_l=l$ by the induction on $l$.
The case $l=0$ is obvious. Suppose that $l\geq 1$ and $\deg\, c_i=i$ 
for $i\leq l-1$. Then the equation (\ref{y-1}) to determine $b_i$, $i\leq l$
can be written as
\bea
&&
1+\sum_{ni+sj<ns}\lambda_{ij}t^{ns-ni-sj}(\sum_{r=0}^{l-1} c_rt^r)^j
=\sum_{i=0}^l b_i t^i
\quad
\hbox{mod. $(t^{l+1})$}.
\label{y-4}
\ena
Since the left hand side of (\ref{y-4}) is of degree $0$ by 
the induction hypothesis, $\deg\, b_i=i$ for any $i\leq l$.
Similarly the equation (\ref{y-3}) determining $c_l$ from $b_1,...,b_l$
can be written as
\bea
&&
\sum_{i=0}^l c_it^i=(\sum_{i=0}^l b_i t^i)^{1/n}
\quad
\hbox{mod. $(t^{l+1})$}
\ena
The right hand side of this equation is of degree $0$. Thus $\deg\, c_l=l$.
\qed
\vskip2mm

\noindent
(ii) This easily follows from (i).
\vskip2mm
\noindent
(iii) By the property (\ref{omega-exp}) for $\homega$,
$$
\homega-\frac{dt_1dt_2}{(t_1-t_2)^2}
$$
is holomorphic near $\{\infty\}\times\{\infty\}$. 
Thus it is possible to expand it as
\bea
&&
\homega(p_1,p_2)-\frac{dt_1dt_2}{(t_1-t_2)^2}
=\left(\sum_{i,j=0}^\infty a_{ij}t_1^it_2^j\right)dt_1dt_2.
\label{homega-exp}
\ena
Let us prove that $a_{ij}$ is a homogeneous polynomial in $\{\lambda_{kl}\}$
of degree $2+i+j$.

Since $\homega(p_1,p_2)(t_1-t_2)^2$ is holomorphic
near $\{\infty\}\times\{\infty\}$ one can write
\bea
&&
\homega(p_1,p_2)(t_1^n-t_2^n)^2
=(\sum_{i,j=0}^\infty b_{ij}t_1^it_2^j)dt_1dt_2.
\non
\ena
We first prove $b_{ij}\in {\mathbb Q}[\{\lambda_{kl}\}]$.

Using 
\bea
&&
f_xdx+f_ydy=0,
\quad
(x_1-x_2)^2=(t_1t_2)^{-2n}(t_1^n-t_2^n)^2,
\non
\ena
we see that
\bea
&&
\homega(p_1,p_2)(t_1^n-t_2^n)^2
=(t_1t_2)^{2n}\,P\, \frac{dx(p_1)dx(p_2)}{f_y(p_1)f_y(p_2)}
\non
\ena
for some homogeneous polynomial $P$ in $x_1,y_1,x_2,y_2, \{\lambda_{kl}\}$ 
with the coefficient in ${\mathbb Q}$ of degree $2s(n-1)$.
By (i), (ii) one can write
\bea
&&
\homega(p_1,p_2)(t_1^n-t_2^n)^2=
\left(\sum_{ij\in {\mathbb Z}}b_{ij}t_1^it_2^j\right)dt_1dt_2,
\quad
b_{ij}\in {\mathbb Q}[\{\lambda_{kl}\}],
\quad
\deg\,b_{ij}=-2n+2+i+j.
\non
\ena
Since the left hand side is holomorphic, the summation in the right hand side
is in fact taken for $i,j\geq 0$.
Therefore one can write
\bea
&&
\left(\homega-\frac{dt_1dt_2}{(t_1-t_2)^2}\right)(t_1^n-t_2^n)^2
=\left(\sum_{i,j=0}^\infty c_{ij}t_1^it_2^j)\right)dt_1dt_2,
\non
\ena
where $c_{ij}\in {\mathbb Q}[\{\lambda_{kl}\}]$ and $\deg\, c_{ij}=-2n+2+i+j$.
By (\ref{homega-exp}) we have
\bea
&&
(t_1^n-t_2^n)^2\sum_{i,j=0}^\infty a_{ij}t_1^it_2^j=\sum_{i,j=0}^\infty c_{ij}
t_1^it_2^j,
\non
\ena
which is equivalent to
\bea
&&
a_{ij}-2a_{i+n,j-n}+a_{i+2n,j-2n}=c_{i+2n,j}.
\label{rec-aij}
\ena
Here we set $a_{ij}=0$ if $i<0$ or $j<0$.

Define $(i,j)<(k,l)$ if and only if $i+j<k+l$ or $i+j=k+l$ and $j<l$.
It defines a total order on the set $\{(i,j)|\,i,j\geq 0\}$.
The equation (\ref{rec-aij}) expresses $a_{ij}$ as a linear combination
of $a_{kl}$ with $(k,l)<(i,j)$, $k+l=i+j$ and $\{c_{rs}\}$ with $r+s=i+j+2n$.
Since $a_{k0}=c_{k+2n,0}$, any $a_{ij}$ is a linear combination of $\{c_{rs}\}$
with $r+s=2n+i+j$.
Thus $a_{ij}$ is in ${\mathbb Q}[\{\lambda_{kl}\}]$ and is homogeneous
of degree $2+i+j$.
\qed

\begin{lemma}\label{te-exp}
(i) The expansion of $\tE(\tilde{p}_1,\tilde{p}_2)$ near 
$(\tinfty,\tinfty)$
is of the form
\bea
&&
\tE(\tilde{p}_1,\tilde{p}_2)=(t_1-t_2)(t_1t_2)^{g-1}
(1+\sum_{i+j\geq 1}c_{ij}t_1^it_2^j),
\non
\ena
where $c_{ij}$ is a homogeneous element in ${\mathbb Q}[\{\lambda_{kl}\}]$
of degree $i+j$.
\vskip2mm
\noindent
(ii) The expansion of $\tE(\tinfty,\tilde{p})$ near $\tinfty$
is of the form
\bea
&&
\tE(\tinfty,\tilde{p})=t^g(1+\sum_{j=1}^\infty c_{0j}t^j),
\non
\ena
where $c_{0j}$ is the same as that in (i).
\end{lemma}
\vskip2mm
\noindent
{\it proof.}
(i) Using the definition, property (iii) of the prime form and (\ref{normal-2})
we have the expansion of the form
\bea
&&
\tE(\tilde{p}_1,\tilde{p}_2)=(t_1-t_2)(t_1t_2)^{g-1}
(1+\sum_{i+j\geq 1}c_{ij}t_1^it_2^j).
\non
\ena
In order prove that $c_{ij}$ has the required properties we use
(\ref{alg-pfn}). The right hand side of (\ref{alg-pfn}) is calculated
in the following way.

Let $\epn=\exp(2\pi i/n)$.
Since $x(p^{(i)})=x(p)=1/t^n$, we take $t^{(i)}=\epn^{-i}t$ as 
a local parameter of $p^{(i)}$ by rearranging $i$ of $p^{(i)}$ if necessary.
Then
\bea
&&
x(p_k^{(i)})=1/t_k^n,
\quad
y(p_k^{(i)})=\frac{\epn^{is}}{t_k^s}(1+\cdots)=
y(p_k)|_{t_k\lar \epn^{-i}t_k}.
\non
\ena
Using these local parameters we get
\vskip5mm
$
\displaystyle{
\exp\left(
\sum_{i=1}^{n-1}
\int_{\tilde{p}_1^{(i)}}^{\tilde{p}_2^{(i)}}\int_{\tilde{p}_1}^{\tilde{p}_2}\homega\right)
}
$
\vskip5mm
$
\displaystyle{
=n^2
\frac{(t_1t_2)^{n-1}}{\prod_{i=1}^{n-1}(t_1-t_2^{(i)})^2}
\exp\left(
\sum_{i=1}^{n-1}
\sum_{k,l=0}^\infty
\frac{a_{kl}}{(k+1)(l+1)}
((t_2^{(i)})^{k+1}-(t_1^{(i)})^{k+1})(t_2^{l+1}-t_1^{l+1})
\right),
}
$
\vskip5mm
and
\bea
&&
\frac{\left(x(p_2)-x(p_1)\right)^2}{f_y(p_1)f_y(p_2)}
=
\frac{1}{n^2}(t_1t_2)^{2g-n-1}(t_1^n-t_2^n)^2
\prod_{j=1}^2\left(1+\sum_{i=1}^\infty c_i't_j^i\right),
\non
\ena
where $c_i'$ is that in Lemma \ref{lem-expansion}, (ii).
The assertions for $c_{ij}$ follow from these expressions and 
Lemma \ref{lem-expansion}.
\qed
\vskip5mm

\noindent
Proof of Theorem \ref{main-cor}.
\vskip2mm
\noindent
(i), (ii): Let $t_i=t(p_i)$. By Lemma \ref{te-exp} we have
\vskip2mm
\bea
&&
\frac{\prod_{i=1}^N\tE(\infty,p_i)^N}{\prod_{i<j}\tE(p_i,p_j)}
=
\frac{(\prod_{i=1}^Nt_i)^{N+g-1}}{\prod_{i<j}(t_i-t_j)}
(1+\sum_{k_1+\cdots+k_N\geq 1}\tilde{c}_{k_1...k_N}t_1^{k_1}\cdots t_N^{k_N}),
\non
\ena
\vskip5mm
\noindent
where $\tilde{c}_{k_1...k_N}\in {\mathbb Q}[\{\lambda_{kl}\}]$
and $\deg\, \tilde{c}_{k_1...k_N}=\sum k_i$.
By (\ref{fi-exp}) we have
\vskip2mm
$
\displaystyle{
(f_1(t),...,f_N(t))
}
$
\vskip5mm
\hskip8mm
$
\displaystyle{
=
(1, \frac{1}{t^{w_2^\ast}}(1+O(t)),...,\frac{1}{t^{w_g^\ast}}(1+O(t)),
\frac{1}{t^{2g}}(1+O(t)),...,\frac{1}{t^{N+g-1}}(1+O(t))),
}
$
\vskip5mm
\noindent
where all $O(t)$ parts are series in $t$ with the coefficients in 
${\mathbb Q}[\{\lambda_{kl}\}]$ and are homogeneous of degree $0$ with
respect to $\{t, \lambda_{kl}\}$.
Using Lemma \ref{gap-seq} we get
\vskip5mm
$
\displaystyle{
(N+g-1,...,N+g-1)+
(0,-w_2^\ast,...,-w_g^\ast,-2g,...,-(N+g-1))
}
$
\vskip5mm
\hskip40mm
$
\displaystyle{
=
(\lambda(n,s)_1,...,\lambda(n,s)_g,0,...,0)
+(N-1,N-2,...,1,0).
}
$
\vskip5mm
\noindent
Let us denote the partition $(\lambda(n,s),0^{N-g})$ by $\lambda^{(N)}(n,s)$.
Then
\bea
&&
\!\!\!\!\!\!\!\!\!\!\!\!
\frac{(\prod_{i=1}^Nt_i)^{N+g-1}}{\prod_{i<j}(t_i-t_j)}
\det(f_i(t_j))_{1\leq i,j\leq N}
=
s_{\lambda^{(N)}(n,s)}(t_1,...,t_N)
+
\sum\widehat{c}_{k_1...k_N}t_1^{k_1}\cdots t_N^{k_N},
\label{Schur-t}
\ena
\vskip2mm
\noindent
where $\widehat{c}_{k_1...k_N}\in {\mathbb Q}[\{\lambda_{kl}\}]$,
$\deg\,\widehat{c}_{k_1...k_N}=-|\lambda(n,s)|+\sum k_i$ and
the summation is taken for $k_i$'s satisfying 
$\sum k_i>|\lambda(n,s)|$.

By (\ref{dui-exp}) we have
\bea
&&
\int_{\infty}^{\tilde{p}} du_{i}=
\frac{t^{w_i}}{w_i}
+\sum_{j=1}^\infty c_{i,j}t^{j+w_i},
\quad
c_{i,j}\in {\mathbb Q}[\{\lambda_{kl}\}],
\quad
\deg\, c_{ij}=j.
\non
\ena
Let 
\bea
&&
T_k=T_k(t_1,...,t_N)=\frac{\sum_{j=1}^Nt_j^k}{k}.
\non
\ena
Then $T_1,...,T_N$ are algebraically independent
and become a generator of the ring of symmetric polynomials of 
$t_1,...,t_N$ with the coefficients in ${\mathbb Q}$,
\bea
&&
{\mathbb Q}[t_1,...,t_N]^{S_N}={\mathbb Q}[T_1,...,T_N].
\non
\ena
Moreover if we prescribe degrees for $t_i$ and $T_i$ by
$$
\deg\, t_i=-1,\qquad \deg\, T_i=-i,
$$ 
a symmetric homogeneous 
polynomial of $t_1,...,t_N$ of degree $k$ can be uniquely written  
as a homogeneous polynomial of $T_1,...,T_N$ of degree $k$.

We have
\bea
u_{i}=\sum_{k=1}^N\int_{\infty}^{\tilde{p}_k} du_i
&=&T_{w_i}+\sum_{j=1}^\infty (w_i+j)c_{ij}T_j
\non
\\
&=&
T_{w_i}+\sum_{\sum jk_j>w_i}\tilde{c}_{k_1...k_N}T_1^{k_1}\cdots
T_N^{k_N},
\ena
where $\tilde{c}_{k_1...k_N}\in {\mathbb Q}[\{\lambda_{kl}\}]$,
 $\deg, \tilde{c}_{k_1...k_N}=-w_i+\sum jk_j$ and
the second expression is unique.

Let us take $N\geq 2g-1=w_g$. Then $T_{w_1}$,...,$T_{w_g}$
are algebraically independent.
We denote the right hand side of (\ref{main-1}) by 
$F(\tilde{p}_1,...,\tilde{p}_N)$.
By Theorem \ref{sigma-1} the series expansion of $F$ in 
$t_1,...,t_N$ can be rewritten as a series in $u_{1},...,u_{g}$.

Let
\bea
&&
F(p_1,...,p_N)=s_{\lambda^{(N)}(n,s)}(t_1,...,t_N)+\sum F_{-k}
\label{decomp-t}
\ena
be the homogeneous decomposition of the series in $t_i$ and
\bea
&&
F(p_1,...,p_N)=\sum_{d\geq 0} \hF_{-d}
\label{decomp-u}
\ena
the homogeneous decomposition of the series in $u_i$, where
in (\ref{decomp-t}) the sum in $k$ is taken for 
$k>|\lambda(n,s)|$.
Let $d_0$ be the smallest integer such that $\hF_{-d}\neq 0$.
Write
\bea
&&
\hF_{-d}(u)=\sum_{\sum l_iw_i=d}\widehat{c}_{d;l_1...l_g}u_{1}^{l_1}
\cdots u_{g}^{l_g}.
\non
\ena
Then
\bea
&&
\hF_{-d_0}(u)=\sum_{\sum l_iw_i=d_0}\widehat{c}_{d_0;l_1...l_g}T_{w_1}^{l_1}
\cdots T_{w_g}^{l_g}+\sum_{\sum jk_j>d_0}c'_{k_1...k_N}T_1^{k_1}\cdots
T_N^{k_N},
\label{F-d0}
\ena
for some constants, with respect to $T_i$'s, $c'_{k_1...k_N}$.
By comparing (\ref{decomp-t}) with (\ref{F-d0}) we have 
$d_0=|\lambda(n,s)|$ and
\bea
&&
\sum_{\sum l_iw_i=d_0}\widehat{c}_{d_0;l_1...l_g}T_{w_1}^{l_1}
\cdots T_{w_g}^{l_g}=
s_{\lambda^{(N)}(n,s)}(t_1,...,t_N).
\non
\ena
Since
\bea
&&
s_{\lambda^{(N)}(n,s)}(t_1,...,t_N)=S_{\lambda(n,s)}(T),
\non
\ena
by (\ref{S-extend}) and (\ref{s=S-1}) and $T_{w_1}$,...,$T_{w_g}$
are algebraically independent, we have
\bea
&&
\hF_{-d_0}(u)=\sum_{\sum l_iw_i=d_0}\widehat{c}_{d_0;l_1...l_g}u_{1}^{l_1}
\cdots u_{g}^{l_g}
=
S_{\lambda(n,s)}(T)|_{T_{w_i}=u_i}.
\non
\ena
For $d>d_0$ one can write 
\bea
&&
F_{-d}(t_1,...,t_N)=
\sum b_{d;l_1...l_N}T_1^{l_1}\cdots T_N^{l_N},
\non
\ena
with $b_{d;l_1...l_N}\in {\mathbb Q}[\{\lambda_{kl}\}]$ and 
$\deg\, b_{d;l_1...l_N}=-d_0+d$.
Comparing this with the expression of $\hF_{-d}(u)$ and
using the algebraic independence of $T_1,...,T_N$,
we see that $\hF_{-d}$ is a polynomial in $\{u_i\}$ with
the coefficients in ${\mathbb Q}[\{\lambda_{kl}\}]$ which
is homogeneous of degree $-d_0$ with respect to $\{u_i,\lambda_{kl}\}$.
Consequently the equation (\ref{main-1}) holds and (i),(ii) 
of Theorem \ref{main-cor}
is proved.
\vskip2mm
\noindent
(iii) Since Riemann's constant $\tau \delta'+\delta''$ 
is a half period, $\sigma(u)$ is even or odd
by Proposition \ref{sigma-anal} and (\ref{R-theta3}). 
Thus the relation in (iii)
 follows from Proposition \ref{Schur-prop2} (iv).
\qed
\vskip5mm

\noindent
{\it Calculation of $C_N$ in Theorem \ref{sigma-2}}
\vskip2mm
We set $t_{1i}=t(\tilde{p}_i)$ and $t_{2i}=t(\tilq_i)$.
We take the limit $\tilq_i\rightarrow \tinfty$, $i=1,...,N$.
Since we know that $F_N$ is holomorphic in $\tilq_i$ by the already proved part
of Theorem \ref{sigma-2}, we can calculate the limit by taking
the limits $t_{21}\rightarrow 0$, $t_{22}\rightarrow0$,... $t_{2N}\rightarrow0$
in this order. Therefore we assume $|t_{21}|<\cdots<|t_{2N}|$ and 
expand everything first into Laurent series 
in $t_{21}$ and take the top term of the expansion. 
Next we expand this top term in $t_{22}$ and so on.
Below, if we simply write $A+\cdots$, then it means such an expansion.
Let us carry out such calculations.

We have 
\bea
&&
D_N=(-1)^{\frac{1}{2}(n-1)N(N-1)}
\epsilon_n^{-\frac{1}{2}N(N-1)-\frac{1}{4}N(N-1)(n-1)^2(n-2)+\frac{1}{2}N(n-1)(2g-n)+\frac{1}{2}gN(n-1)(n-2)}
\non
\\
&&
\qquad\quad
\times 
\epsilon(1)^N\det(f_i(p_j))_{1\leq i,j\leq N}
\prod_{j=1}^N t_{2j}^{(j-1)(n-1)^2-\sum_{i=1}^{n-1}(nN+g-i)}
+\cdots.
\label{s2-1}
\ena
By Lemma \ref{te-exp} (i) we have
\bea
&&
\!\!\!\!\!\!\!\!\!\!\!\!
M_N=(-1)^N
\frac{
\prod_{i=1}^N(t_{1i}t_{2i})^{nN+g-1}\prod_{i,j=1}^N(t_{1i}-t_{2j})
}
{
\prod_{i<j}(t_{1i}-t_{1j})(t_{2i}-t_{2j})
\cdot \prod_{i,j=1}^N(t_{1i}^n-t_{2j}^n)}
\left(1+(\hbox{regular term})\right).
\label{s2-2}
\ena
\vskip5mm
By calculation
\vskip5mm
\hskip5mm
$
\displaystyle{
F_N/D_N=(-1)^{\frac{1}{2}N(N-1)(n-2)}\epsilon(s)^{-N}
}
$
\bea
&&
\!\!\!\!
\times
\frac{(\prod_{i=1}^Nt_{2i})^{(N-1)n(n-2)+\frac{1}{2}(n-1)(n-2)s}}
{\prod_{i<j}(t_{2i}^n-t_{2j}^n)^{n-2}}
\left(1+(\hbox{regular term})\right).
\label{s2-3}
\ena
\vskip2mm
Let $H_N$ be the product of $F_N/D_N$ and $M_N$. Then
\vskip5mm
\hskip10mm
$
\displaystyle{
H_N=(-1)^{\frac{1}{2}N(N+1)}\epsilon(s)^{-N}
\frac{(\prod_{i=1}^Nt_{1i})^{N+g-1}}
{\prod_{i<j}(t_{1i}-t_{tj})}
}
$
\bea
&&
\times
\prod_{i=1}^Nt_{2i}^{-(i-1)(n-1)^2+nN+g-1+(N-1)n(n-2)+\frac{1}{2}(n-1)(n-2)s}
+\cdots.
\label{s2-4}
\ena
Multiplying (\ref{s2-1}) and (\ref{s2-4}) we have
\bea
M_N F_N&=&C_N^{-1}\,\,\frac{(\prod_{i=1}^Nt_{1i})^{N+g-1}}
{\prod_{i<j}(t_{1i}-t_{tj})}\,\,\det(f_i(p_j))_{1\leq i,j\leq N}
+\cdots.
\non
\ena
Thus, by (\ref{Schur-t}), in the limit $t_{21}\rightarrow0$,...,
$t_{2N}\rightarrow0$ we get
\bea
M_NF_N&=&
C_N^{-1}\,\,s_{\lambda^{(N)}(n,s)}(t_1,...,t_N)+\cdots
\non
\\
&=&C_N^{-1}\,\,S_{\lambda(n,s)}(T)|_{T_{w_i}=u_i}+\cdots,
\non
\ena
where in the last $+\cdots$ part is a series in $t_{1i}$, $i=1,...,N$.
\qed

\subsection{Example$-g=1-$}

We take 
$f(x,y)=y^2-x^3-\lambda_{10}x-\lambda_{00}$, that is, we set 
$\lambda_{2,0}=\lambda_{01}=\lambda_{11}=0$. 
In this case  $c_{ij;kl}$ satisfying (ii) of Proposition \ref{cijkl} is unique,
$c_{00;10}=1$ and other $c_{ij;kl}$'s are all zero. Then
\bea
\homega(p_1,p_2)&=&
d_{p_2}\left(\frac{y_1+y_2}{2y_1(x_1-x_2)}dx_1\right)+\frac{dx_1}{2y_1}
\frac{x_2dx_2}{2y_2}
\non
\\
&=&
\frac{2y_1y_2+x_1x_2(x_1+x_2)+\lambda_{10}(x_1+x_2)+2\lambda_{00}}
{4y_1y_2(x_1-x_2)^2}dx_1dx_2,
\non
\ena
and
\bea
&&
du_1=-\frac{dx}{2y},
\quad
dr_1=-\frac{xdx}{2y}.
\non
\ena
Theorem \ref{sigma-2}
gives
\bea
&&
\sigma\left(\int_{\tilde{q}_1}^{\tilde{p}_1} du_1\right)
=\tE(\tilde{p}_1,\tilde{q}_1)=
\frac{x(q_1)-x(p_1)}{2\sqrt{y(p_1)y(q_1)}}
\exp\left(
\frac{1}{2}
\int_{\tilde{q}_1^{(1)}}^{\tilde{p}_1^{(1)}}
\int_{\tilde{q}_1}^{\tilde{p}_1}\homega
\right),
\non
\ena
where in the right hand side the square root and the value of
the exponential part should be appropriately specified.
Notice that $p^{(1)}=(x,-y)$ for $p=(x,y)$.
One can transform the defining equation of the elliptic curve to 
Weierstrass form
\bea
Y^2=4X^3-g_2X-g_3,
\non
\ena
by
\bea
&&
X=x,
\quad
Y=-2y,
\quad
g_2=-4\lambda_{10},
\quad
g_3=-4\lambda_{00}.
\non
\ena
The sigma function in this case coincides with that of Weierstrass.
Let $\wp(u)$ be the Weierstrass elliptic function.
The symplectic basis $du_1$, $dr_1$ are transformed to
\bea
&&
du=\frac{dX}{Y},
\quad
dr=\frac{XdX}{Y}.
\non
\ena
and $\homega$ becomes
\bea
\homega&=&\frac{2Y_1Y_2+4X_1X_2(X_1+X_2)-g_2(X_1+X_2)-2g_3}
{4Y_1Y_2(X_1-X_2)^2}dX_1dX_2
\non
\\
&=& 
\wp(v_2-v_1)dv_1dv_2,
\non
\ena
where 
$$
v_i=\int_{\tinfty}^{\tilde{p}_i}du.
$$
The formula for the sigma function becomes
\bea
&&
\sigma\left(v_2-v_1\right)=
\frac{X_1-X_2}{\sqrt{Y_1Y_2}}
\exp\left(
\frac{1}{2}
\int_{\tilde{p}_1^{(1)}}^{\tilde{p}_2^{(1)}}\int_{\tilde{p}_1}^{\tilde{p}_2}\homega
\right).
\non
\ena
This formula coincides with that given by Klein \cite{K1}.

\subsection{Example-Hyperelliptic Case-}
Consider the case
\bea
&&
f(x,y)=y^2-x^{2g+1}-\sum_{i=0}^{2g}\lambda_{i,0}x^i.
\non
\ena
We set $\lambda_{2g+1,0}=1$.
One can take $c_{i_1,0;i_2,0}=(i_2-i_1)\lambda_{i_1+i_2+2,0}$
for $0\leq i_1\leq g-1$, $i_1+1\leq i_2\leq 2g-i_1$ 
and all other $c_{i_1j_1;i_2j_2}$ to be zero.
Then \cite{BEL1,B2}
\bea
\homega&=&
d_{p_2}\left(\frac{y_1+y_2}{2y_1(x_1-x_2)}\right)
+\sum_{i=1}^g du_i(p_1)dr_i(p_2)
\non
\\
&=&\frac{
2y_1y_2+
\sum_{i=0}^gx_1^ix_2^i\left(2\lambda_{2i,0}+\lambda_{2i+1,0}(x_1+x_2)
\right)}
{
4y_1y_2(x_1-x_2)^2
}dx_1dx_2,
\non
\ena
and
\bea
du_i&=&-\frac{x^{g-i}dx}{2y},
\non
\\
dr_i&=&
-
\sum_{i=g+1-i}^{g+i}(k-g+i)\lambda_{k+g+2-i,0}
\frac{x^kdx}{2y}.
\non
\ena
Set $p_i=(x_i,y_i)$ and $q_i=(X_i,Y_i)$.
The formula of the sigma function is given by (\ref{general}) with
\bea
&&
F_N=
\left|
\begin{array}{cccccccccc}
1&\cdots&1&1&\cdots&1&\\
x_1&\cdots&x_N&X_1&\cdots&X_N&\\
\vdots&\quad&\vdots&\vdots&\quad&\vdots&\\
x_1^g&\cdots&x_N^g&X_1^g&\cdots&X_N^g&\\
y_1&\cdots&y_N&-Y_1&\cdots&-Y_N&\\
x_1^{g+1}&\cdots&x_N^{g+1}&X_1^{g+1}&\cdots&X_N^{g+1}&\\
x_1y_1&\cdots&x_Ny_N&-X_1Y_1&\cdots&-X_NY_N&\\
\vdots&\quad&\vdots&\vdots&\quad&\vdots&\\
\end{array}
\right|,
\non
\ena
and
\bea
&&
C_N=(-1)^{\frac{1}{2}N(N+1)+gN}.
\non
\ena
This is Klein's formula \cite{K1,K2}.

\section{Concluding Remarks}
In this paper we have established the formula for
the sigma function associated to a $(n,s)$-curve $X$ in terms of algebraic
integrals. Some properties of the series expansion of the sigma function
are deduced from it. Namely it is shown that the first term of the expansion
becomes the Schur function corresponding to the partition determined from
the gap sequence at infinity and the expansion coefficients are homogeneous
polynomials of the coefficients of the defining equation of the curve. 

The building block of the formula is the prime function.
It is a multi-valued function on $X\times X$ with some vanishing property
and has the same transformation
rule as that of the sigma function if one of the variables goes round 
a cycle of $X$. Remarkably,
in the case of hyperelliptic curves, \^Onishi \cite{O1} has constructed 
a function with the same properties as a certain derivative of the 
sigma function. By the uniqueness of such a function they coincide.
In general it is expected that the prime function can be expressed 
as a derivative of the sigma function. 

Fay's determinant formula ((43) in \cite{F}) expresses Riemann's theta
function in terms of prime form and some determinant. To get a formula
of the sigma function one needs to take a limit sending some parameter
to a singular point of the theta divisor. To this end one has to know
the structure of sigularities of the theta divisor. So it is difficult
to get the formula of the sigma function by taking a limit of Fay's
formula in general. We have avoided this task and directly constructed 
a formula of the sigma function. 
 
\vskip10mm

\noindent
{\large {\bf Acknowledgement}} 
\vskip3mm
\noindent
The author would like to thank Victor Enolski
for answering questions and valuable comments and to Koji Cho, Yasuhiko
Yamada for useful discussions. The author is also grateful to 
Shigeki Matsutani and Yoshihiro \^Onishi for their encouragement.
This research is supported by Grant-in-Aid for Scientific Research (B) 
17340048.

\end{document}